\documentclass[10pt]{amsart}

\usepackage{amsmath}
\usepackage{amssymb}
\usepackage[all]{xy}

\makeatletter
\@addtoreset{equation}{section}
\makeatother

\newtheorem{theorem}{Theorem}[section]
\newtheorem{proposition}[theorem]{Proposition}

\newtheorem{corollary}[theorem]{Corollary}
\newtheorem{conjecture}[theorem]{Conjecture}

\theoremstyle{definition}
\newtheorem{example}[theorem]{Example}

\def\fA{\mathcal{A}}
\def\fB{\mathcal{B}}
\def\fC{\mathcal{C}}
\def\fE{\mathcal{E}}
\def\fF{\mathcal{F}}
\def\fH{\mathcal{H}}

\def\fL{\mathcal{L}}
\def\fP{\mathcal{P}}

\def\fU{\mathcal{U}}
\def\fW{\mathcal{W}}
\def\fX{\mathcal{X}}

\def\Cee{\mathbb{C}}
\def\Dee{\mathbb{D}}
\def\Hee{\mathbb{H}}
\def\En{\mathbb{N}}

\def\Que{\mathbb{Q}}
\def\Ree{\mathbb{R}}
\def\Tee{\mathbb{T}}
\def\Zee{\mathbb{Z}}

\def\O{\mathrm{O}}
\def\Sp{\mathrm{Sp}}
\def\U{\mathrm{U}}
\def\Re{\mathrm{Re}}
\def\Im{\mathrm{Im}}

\def\gg{\mathfrak{g}}

\def\alp{\alpha}
\def\del{\delta}
\def\eps{\varepsilon}
\def\gam{\gamma}
\def\lam{\lambda}
\def\ome{\omega}
\def\sig{\sigma}
\def\Sig{\Sigma}
\def\vphi{\varphi}

\def\aand{\text{ and }}
\def\iif{\text{ if }}
\def\iin{\text{ in }}
\def\ffor{\text{ for }}

\def\iff{\Leftrightarrow}

\def\setdif{\setminus}

\def\endpf{{\hfill$\square$\medskip}}
\def\proof{{\noindent{\bf Proof.}\thickspace}}

\def\mult{\!\cdot\!}
\def\comp{\raisebox{.2ex}{${\scriptstyle\circ}$}}
\def\con{\!*\!}
\def\cross{\!\times\!}
\def\spn{\mathrm{span}}
\def\lap{\raisebox{.2ex}{${\scriptstyle\square}$}}
\def\rap{\raisebox{.2ex}{${\scriptstyle\lozenge}$}}
\def\til#1{\tilde{#1}}

\def\wbar#1{\overline{#1}}
\def\what#1{\widehat{#1}}
\def\inprod#1#2{\left\langle #1 | #2 \right\rangle}
\def\id{\mathrm{id}}
\def\Re{\mathrm{Re}}
\def\Im{\mathrm{Im}}
\def\spn{\mathrm{span}}

\def\alg{\mathrm{alg}}
\def\bl{\mathrm{L}}
\def\ecoef{\mathrm{E}}
\def\ccoef{\mathrm{A}}
\def\wcoef{\mathrm{B}}
\def\falg{\mathrm{A}(G)}
\def\fsal#1{\mathrm{B}(#1)}
\def\fsalg{\mathrm{B}(G)}
\def\foalg{\mathrm{B}_0(G)}
\def\ideal{\mathrm{i}}
\def\meas{\mathrm{M}}

\def\cbmg{\mathrm{M}_{cb}\mathrm{A}(G)}

\def\cstarg{\mathrm{C}^*(G)}
\def\wstarg{\mathrm{W}^*(G)}
\def\vn{\mathrm{VN}}
\def\opalg{\mathrm{OA}}
\def\posd{\mathrm{P}}
\def\pol{\mathrm{pol}}
\def\queg{\mathrm{Q}(G)}
\def\wc{\leq_w}
\def\qc{\leq_q}
\def\mq{\leq_{mq}}
\def\cmq{\leq_{\bar{m}q}}
\def\ec{\leq_e}
\def\me{\leq_{me}}
\def\se{\leq_{se}}
\def\we{\leq_{we}}

\def\unorm#1{\left\|#1\right\|_\infty}
\def\fsnorm#1{\left\|#1\right\|_{\mathrm B}}
\def\norm#1{\left\|#1\right\|}

\begin{document}

\title[Matrix coefficients and compactifications]
{Matrix coefficients of unitary representations and associated compactifications}

\author{Nico Spronk and Ross Stokke}

\begin{abstract}
We study, for a locally compact group $G$, the compactifications $(\pi,G^\pi)$
associated  with unitary representations $\pi$, which we call
{\it $\pi$-Eberlein compactifications}.  We also study the Gelfand spectra $\Phi_{\fA(\pi)}$
of the uniformly closed algebras $\fA(\pi)$ generated by matrix coefficients of such $\pi$.
We note that $\Phi_{\fA(\pi)}\cup\{0\}$ is itself a semigroup and
show that the \v{S}ilov boundary of $\fA(\pi)$ is $G^\pi$.  We study containment
relations of various uniformly closed algebras generated by matrix coefficients, and
give a new characterisation of amenability:  the constant function $1$ can be uniformly
approximated by matrix coefficients of representations weakly contained in
the left regular representation if and only if $G$ is amenable.
We show that for the universal
representation $\ome$, the compactification $(\ome,G^\ome)$ has a certain universality property:
it is universal amongst all compactifications of $G$ which may be embedded as
contractions on a Hilbert space, a fact which was also recently proved by
Megrelishvili \cite{megrelishvili}.  We illustrate our results with examples including various abelian and
compact groups, and the $ax+b$-group.  In particular, we witness algebras
$\fA(\pi)$, for certain non-self-conjugate $\pi$, as being generalised algebras of
analytic functions.
\end{abstract}

\maketitle

\footnote{{\it Date}: \today.

2000 {\it Mathematics Subject Classification.} Primary 43A30, 47D03, 46J10,
43A07;
Secondary 43A25, 43A65, 46E25.
{\it Key words and phrases.} Fourier-Stieltjes algebras, semitopological compactification,
\v{s}ilov boundary, amenable group.

The first-named author is supported by NSERC Grant 312515-2010. The second named
author is supported by NSERC Grant 298444-2010.}


\section{Preliminaries}

\subsection{Introduction}
Given a locally compact group $G$, introverted subspaces of $\fC\fB(G)$, the
continuous bounded functions on $G$,
and their associated compactifications of $G$ have been studied by many
authors over the years; see, for example,
the treatise of Berglund, Junghenn and Milnes \cite{berglundjm}.
Let $\fsalg$ denote the Fourier-Stieltjes algebra of $G$ and $\fE(G)$
its uniform closure in $\fC\fB(G)$, which we call the {\it Eberlein algebra}.
The algebras of left uniformly continuous functions
$\fL\fU\fC(G)$, weakly almost periodic functions $\fW\fA\fP(G)$
and almost periodic functions
$\fA\fP(G)$, and their associated compactifications $G^{\fL\fU\fC}$,
$G^{\fW\fA\fP}$ and $G^{\fA\fP}$, have
received a great deal of attention over the years, while less has been
paid to $\fE(G)$ and $G^\fE$.  Being a quotient of $\fW\fA\fP(G)^*$, $\fE(G)^*$
inherits many of the nicest properties of $\fW\fA\fP(G)^*$, in particular Arens
regularity.  In some situations $\fE(G)=\fW\fA\fP(G)$; for connected
$G$ this has been characterised by Mayer \cite{mayer}.

In the present article we initiate a comprehensive investigation
into the properties of $\fE(G)$, $\fE(G)^*$ and the associated compactification
$G^\fE$.  This is a natural task as
Eberelein \cite{eberlein} initiated the theory of weakly almost periodic functions
on abelian groups in order to gain understanding of the
Fourier-Stieltjes transforms of measures on the dual group.
We note that the spectrum $\Phi_{\fsalg}$ of $\fsalg$
can be very complicated, especially for abelian groups --- see the monographs
of Rudin \cite{rudinB} and Graham and McGehee \cite{grahamm}.
Meanwhile, the spectrum of $\fE(G)$, being the closure of $G$ in
$\Phi_{\fsalg}$, is much more tractable.
In order to be systematic, we restrict ourselves not only to examining
$\fE(G)$, but, in fact uniform algebras $\fA(\pi)$ generated by matrix coefficients
of general unitary representations $\pi$.  Herein we gain some extra
complications and generalise, in a certain manner, the theory of
``generalised algebras of analytic functions" in the sense of Arens and Singer
\cite{arenss}.  We find, for example, that the \v{S}ilov boundary of
$\fA(\pi)$ is $G^{\fA(\pi)}\setdif\{0\}$, where $G^{\fA(\pi)}$
is the associated compactification of $G$.

One of the most surprising features is that we can use semigroup
structure theory to gain an analogue of the Fell-Hulanicki characterisation
of amenability.  We also show that $G^\fE$ is the universal compactification
amongst those compactifications of $G$ which are realisable as
contractions on a Hilbert space, a result already proved by Megrelishvili
\cite{megrelishvili}, although with different techniques.
An interesting feature which arises is that there are involutive
compactifications of $G$, i.e.\ admiting an involution $x\mapsto x^*$
which extends $s\mapsto s^{-1}$, which cannot faithfully represented
as contractions on a Hilbert space.

The second named author owes his interest in Eberlein compactifications
to his work on group algebra homomorphism problems.
In \cite{Stokke}, it is established that  $*$-homomorphisms of 
$L^1(G)$ into  $M(H)$ are in bijective correspondence with
weak$^*$-continuous $*$-homomorphisms from $M(G^\fE)$ into $M(H^\fE)$, and 
from $M(G^\fE)$ into $M(H)$.

\subsection{Plan}

While our focus and goal is to understand these certain function spaces
for locally compact groups
we have realised that much of the general theory can be framed in the much
more general context of a semitopological semigroup.  Hence in Section
\ref{sec:semigroups} we study spaces of functions over a semi-topological
semigroup, which we need not even assume is locally compact.  The philosophy
of this section is that of \cite{berglundjm}, in which the duality between
certain translation-invariant unital C*-subalgebras of
functions and compactifications is the major tool.  We augment this in
a modest but critical manner.  Since our goal is to understand
compactifications associated to matrix coefficients, as introduced in
\S\S\ref{ssec:basicspaces}, we find it handy to use the notion of
homogeneous spaces of functions on a semigroup, which we describe
in \S\S\ref{ssec:homogeneousspaces}.  We give a systematic exposition
of the basic theory of such subspaces and discuss the algebras and
self-adjoint algebras generated by them.  In \S\S\ref{ssec:semitopological}
we focus on semitopological compactifications.  Particularly, we observe
that if the underlying semigroup $G$ itself has a continuous involution ---
say $s\mapsto s^{-1}$ in the case of a topological group --- then the
weakly almost periodic compactification is universal amongst
compactifications which themselves admit a continuous involution
extending that of $G$.  In \S\S\ref{ssec:CH} we introduce
the concept of  (CH)-compactifications, those realisable as weak*-closed
semigroups of contractions on a Hilbert space.

The heart of this article is Section \ref{sec:eberlein}.  We return to
the setting of a locally compact group $G$.  We study the relationships
between various spaces generated by matrix coefficients, both uniformly
closed and closed in the Fourier-Stieltjes norm.  We introduce the
concept of {\it Eberlein containment} and illustrate its relationship
to more classical methods of comparing unitary representations.
In particular we observe that the Eberlein compactification, which is the
spectrum of the uniform closure of the Fourier-Stieltjes algebra in $\fC\fB(G)$,
is an invariant for $G$.  In \S\S\ref{ssec:subalggenmatcoef}
we study the spectra of algebras generated by matrix coefficients,
giving a criterion for determining elements of the spectra of these algebras
in the vein
of Walter \cite{walter,walter1}, and characterising the \v{S}ilov boundary
within the spectrum.  In \S\S\ref{ssec:ebweakcontamen} we
manipulate the role played by almost periodic
functions to characterise amenability of $G$ in terms of some
uniformly closed algebras of functions.  In \S\S\ref{ssec:universal}
we recover a result already shown by Megrelishvili \cite{megrelishvili},
that describes a natural universality property of the Eberlein compactification.
As it is further noted in \cite{megrelishvili}, there exist monothetic
compact semitopological semigroups which are not (CH)-compactifications.
We extend this observation to include wider classes of groups, using
results of Chou \cite{chou} and Mayer \cite{mayer,mayer1}.

In Section \ref{sec:examples} we illustrate aspects of our theory
with examples.  We include spine-type examples, after
Ilie and the first named author \cite{ilies} and Berglund \cite{berglund},
however we modify them to show special properties of the compactifications.
We show how abelian groups fit into our theory and even compute
spectra for subsemigroups of integers and open subsemigroups in vector groups.
This emphasises that the non-self-adjoint algebras of matrix coefficients
are indeed ``generalised algebras of analytic functions'' in the sense of \cite{arenss}.
We continue on this track by illustrating
an example for compact groups, and finally the $ax+b$-group.

\subsection{The basic spaces}\label{ssec:basicspaces}
We consider several spaces of functions based upon unitary representations
of a locally compact group $G$.  We let $\Sigma_G$ denote the class of
all continuous unitary representations $\pi$, continuous in the
sense that each matrix coefficient function, $s\mapsto\inprod{\pi(s)\xi}{\eta}$,
is continuous on $G$.  For two elements $\pi$, $\sig$ of
$\Sig_G$, we write $\pi\cong\sig$ to denote the relation of unitary equivalence.
We let $\{\bar{\pi},\bar{\fH}_\pi\}$ denote the conjugate representation.
As in \cite{arsac}, for $\pi\iin\Sig_G$ we define
\[
F_\pi=\spn\{\inprod{\pi(\cdot)\xi}{\eta}:\xi,\eta\in\fH_\pi\}.
\]
We observe that $F_\pi$ is clearly left and right translation invariant and that
\begin{equation}\label{eq:inversioninvariance}
F_\pi^\vee=\bar{F}_\pi=F_{\bar{\pi}}
\end{equation}
where $\check{u}(s)=u(s^{-1})$.  Indeed
$\inprod{\pi(\cdot)\xi}{\eta}^\vee=\wbar{\inprod{\pi(\cdot)\eta}{\xi}}=
\inprod{\bar{\pi}(\cdot)\bar{\eta}}{\bar{\xi}}$
for $\xi,\eta\iin\fH_\pi$.  Hence $F_\pi$ is inversion-invariant if $\pi\cong\bar{\pi}$.

We let $\fsalg=\bigcup_{\pi\in\Sig_G} F_\pi$ denote the {\it Fourier-Stieltjes
algebra} of $G$ as defined in \cite{eymard}.
This is an algebra of functions and, moreover a Banach
algebra when endowed with the norm admitting the two equivalent descriptions
below
\begin{align*}
\fsnorm{u}&=\inf\bigl\{\norm{\xi}\norm{\eta}:u=\inprod{\pi(\cdot)\xi}{\eta},
\pi\in\Sig_G,\xi,\eta\in\fH_\pi\bigr\} \\
&=\sup\left\{\left|\int_G u(s)f(s)ds\right|:f\in\bl^1(G),\norm{f}_*=
\sup_{\pi\in\Sig_G}\norm{\pi_1(f)}\leq 1\right\}
\end{align*}
where $\pi_1(f)=\int_Gf(s)\pi(s)ds$ (weak* integral).
We note that we obtain the duality identification
$\fsalg\cong\cstarg^*$, where $\cstarg=\wbar{\bl^1(G)}^{\norm{\cdot}_*}$,
the completion of $\bl^1(G)$ in the largest C*-norm $\norm{\cdot}_*$.  
Let $\posd_1(G)$ denote the set of continuous, positive definite functions $u$
with $u(e)=1$.
Then the Gelfand-Naimark construction provides for each $u$ in $\posd_1(G)$
a unitary representation $\pi_u$ and a norm one cyclic vector $\xi$ in $\fH_{\pi_u}$
for which $u=\inprod{\pi_u(\cdot)\xi}{\xi}$.  Moreover, every such cyclic
representation $\{\pi,\xi\}$ arises, up to unitary equivalence, in this manner.
See \cite[(3,20)]{folland} or \cite[13.4.5]{dixmier} for details.
We define the {\it universal representation}
\[
\ome_G=\bigoplus_{u\in\posd_1(G)}\pi_u.
\]
As shown in \cite{eymard}, $\fsalg=\spn\posd_1(G)=F_{\ome_G}$.  Also, it is well-known that
$\norm{(\ome_G)_1(f)}=\norm{f}_*$ for $f\iin\bl^1(G)$.

We fix an element $\pi$ of $\Sig_G$.  We let
\[
\ccoef_\pi=\wbar{F_\pi}^{\fsnorm{\cdot}}\quad\text{and}\quad
\fE_\pi=\wbar{F_\pi}^{\unorm{\cdot}}
\]
which are each closed translation-invariant subspaces of $\fsalg$ and $\fC\fB(G)$, respectively.
We then let $\alg(F_\pi)$ denote the algebra of functions generated by $F_\pi$ and define
\[
\ccoef(\pi)=\wbar{\alg(F_\pi)}^{\fsnorm{\cdot}}\quad\text{and}\quad
\fA(\pi)=\wbar{\alg(F_\pi)}^{\unorm{\cdot}}.
\]
which are translation-invariant closed subalgebras of $\fsalg$ and $\fC\fB(G)$, respectively.
Finally, we let
\begin{gather*}
\ecoef(\pi)=\wbar{\alg(F_\pi+\bar{F}_\pi)}^{\fsnorm{\cdot}},\quad
\fE(\pi)=\wbar{\alg(F_\pi+\bar{F}_\pi)}^{\unorm{\cdot}} \\
\text{and}\quad
\fE_1(\pi)=\wbar{\alg(\Cee 1+F_\pi+\bar{F}_\pi)}^{\unorm{\cdot}}
\end{gather*}
which are  translation-invariant, conjugation-invariant closed subalgebras of $\fsalg$ and 
$\fC\fB(G)$, respectively.  We define representations
\[
\tau_\pi=\bigoplus_{n\in\En}\pi^{n\otimes}\quad\text{and}\quad
\rho_\pi=\bigoplus_{\substack{m,n\in\{0\}\cup\En \\ m+n\geq 1}}
\pi^{m\otimes}\otimes\bar{\pi}^{n\otimes}
\]
where $\tau_\pi$ is defined on the Hilbertian direct sum
$\fH_{\tau_\pi}=\ell^2\text{-}\bigoplus_{n\in\En}\fH_\pi^{n\otimes_2}$, and
$\fH_{\rho_\pi}$ is defined similarly.  Then we have that
\[
\ccoef(\pi)=\ccoef_{\tau_\pi}\quad\text{and}\quad
\ecoef(\pi)=\ccoef_{\rho_\pi}.
\]
For details see \cite[Lem.\ 4.1]{stokke0}.  The fact that $\fsnorm{\cdot}\geq\unorm{\cdot}$
then implies that
\[
\fA(\pi)=\fE_{\tau_\pi},\quad\fE(\pi)=\fE_{\rho_\pi}\quad\text{and}\quad
\fE_1(\pi)=\fE_{1\oplus\rho_\pi}.
\]
We let $\wcoef_\pi$ denote the weak* closure of $F_\pi$ in $\fsalg$.
We observe that it follows from \cite[(1.20)]{eymard} that the representation
$\ome_\pi=\bigoplus_{u\in\posd_1(G)\cap\wcoef_\pi}\pi_u$ satisfies $\wcoef_\pi=F_{\ome_\pi}$.
We define the {\em weak $\pi$-Eberlein algebra} by
\[
\fE\fB(\pi)=\wbar{\alg(\wcoef_\pi+\bar{\wcoef}_\pi)}^{\unorm{\cdot}}=\fE(\ome_\pi).
\]
The definition of the weak $\pi$-Eberlein algebra
is arguably the least natural one here as it mixes topologies; it
is motivated by its use in Theorem \ref{theo:amenrep}.
Finally, we let the {\em Eberlein algebra} of $G$ be given by
\[
\fE(G)=\fE(\ome_G)=\wbar{\fsalg}^{\unorm{\cdot}}
\]
where $\ome_G$ is the universal representation, defined above.

\section{Function spaces over semigroups and compactifications}
\label{sec:semigroups}

For this section we will always let $G$ be a semitopological semigroup,
not necessarily locally compact.

\subsection{Arens products and semigroup compactifications}\label{sses:arensprod}
Our standard reference for this section is
the text \cite{berglundjm}, though our notation differs slightly.

A {\it right topological
compactification} of $G$ is a pair $(\del,S)$, where $S$ is a compact right
topological semigroup --- for any $t\iin S$, $s\mapsto st$ is
continuous --- and $\del:G\to S$ is a continuous homomorphism
whose range is both dense in $S$ and contained in the topological centre
$Z_T(S)=\{t\in G:s\mapsto ts\text{ is continuous}\}$.
We define left topological compactifications similarly.
If $S$ is semitopological, i.e.\ $S=Z_T(S)$, then we say $(\del,S)$
is a {\it semitopological compactification} of $G$.

If $(\del,S),(\eps,T)$ are two right [left] topological compactifications of $G$ we
write $(\del,S)\leq(\eps,T)$ if there is a continuous homomorphism $\theta:
T\to S$ such that $\theta\comp\eps=\del$; necessarily, $\theta$ is surjective.
We say $(\del,S)$ is a {\it factor}
of $(\eps,T)$, conversely $(\eps,T)$ is an {\it extension} of $(\del,S)$.
We say $(\del,S)$ and $(\eps,T)$ are equivalent, written $(\del,S)\cong(\eps,T)$
if $\theta$, above, is an isomorphism.  This condition is the same as simultaneously having
$(\del,S)\leq(\eps,T)$ and $(\del,S)\geq(\eps,T)$.  In particular, $\leq$ is a partial
ordering on the class, in fact the set (see the remark after Theorem \ref{theo:comparison}, below)
of equivalence classes of right [left] topological compactifications of $G$.

Let $\fC\fB(G)$ denote the C*-algebra of continuous complex-valued bounded functions
on $G$ with uniform norm $\unorm{\cdot}$.
If $f\in\fC\fB(G)$, $s\in G$ we denote the
anti-action of left translation and the action of right translation of $s$ on $f$ by
\[
f\mult s(t)=f(st) \aand s\mult f(t)=f(ts)
\]
for $t\iin G$.  Let $\fX$ be a closed linear subspace of $\fC\fB(G)$.
We say $\fX$ is {\it left [right] introverted} if it is closed under left [right] translations,
and for any $m\iin\fX^*$, $m\mult f$, defined by
\[
m\mult f(s)=m(f\mult s)\quad [f\mult m(s)=m(s\mult f)]
\]
is also an element of $\fX$.   Note that we do not insist that $\fX$
contains the constant functions.   We say that
$\fX$ is {\it introverted} if it is both left and right introverted.

If $\fX$ is left [right] introverted then we define the {\it left [right] Arens products}
on $\fX^*$ by
\begin{equation}\label{eq:arensprod}
m\lap n(f)=m(n\mult f)\quad [m\rap n(f)=n(f\mult m)]\quad\ffor f\iin\fX
\end{equation}
which makes $\fX^*$ into a right [left] dual Banach algebra
in the sense that for a fixed $n$, $m\mapsto m\lap n$ [$m\mapsto n\rap m$]
is weak*-weak* continuous on $\fX^*$.
If $\fX$ is introverted, we say that $\fX$ is {\it Arens
regular} if the left and right Arens products coincide.  Arens regularity
is discussed in greater detail in Section
\ref{ssec:semitopological}.

If $\fX$ is a closed subalgebra of $\fC\fB(G)$, then
we denote its Gelfand spectrum by $\Phi_\fX$, and endow it with its weak* topology.
We let $\eps_\fX:G\to\Phi_\fX$ denote the evaluation map, which has dense range
in the case that $\fX$ is a C*-algebra.
We say that $\fX$ is {\it left [right] m-introverted} if
$\chi\mult f\in\fX$ [$f\mult\chi\in\fX$] for each $\chi\iin\Phi_\fX$.
We record, for ease of reference, the following standard result which can be
found as Theorems 3.1.7 and 3.1.9 in \cite{berglundjm}.

\begin{theorem}\label{theo:comparison}
{\bf (i)} Let $\fX$ be a left [right] translation invariant unital C*-subalgebra
of $\fC\fB(G)$. Then $\fX$ is left [right] m-introverted
if and only if $\Phi_\fX$ is a right [left] topological semigroup under the
left [right] Arens product of (\ref{eq:arensprod}).  In this case
$(\eps_\fX,\Phi_\fX)$ is a right [left] topological
compactifiaction of $G$ and  $\fX=\fC(\Phi_\fX)\comp\eps_\fX$.

{\bf (ii)} If $(\del,S),(\eps,T)$ are two right [left] topological compactifications
of $G$ then
\[
(\del,S)\leq(\eps,T)\quad\text{ if and only if }\quad\fC(S)\comp\del\subset\fC(T)\comp\eps.
\]
\end{theorem}

In particular, the family of all equivalence classes of
right [left] topological compactifications of $G$
is a set, realised in bijective correspondence with the left [right] m-introverted
unital C*-subalgebras of $\fC\fB(G)$.




\subsection{Homogeneous subspaces and the algebras they generate}\label{ssec:homogeneousspaces}
Let us consider a mild generalisation of the concept of introverted subspaces
of $\fC\fB(G)$, which will be useful for our goals.
For the sake of brevity we will work mainly with left actions on subspaces and
hence right topological dual algebras and compactifications; the opposite handed
analogues are similar.
A {\it left homogeneous subspace} of $\fC\fB(G)$ is a subspace
$X$ such that
\begin{itemize}
\item[(i)] $X$ is equipped with a norm $\norm{\cdot}$ under which it is complete and
for which $\norm{f}\geq\unorm{f}$ for $f\iin X$; and
\item[(ii)] $X$ is left translation invariant and
$(f,s)\mapsto f\mult s:X\cross G\to X$ is continuous in $s$ and
contractive in $f$.
\end{itemize}
Moreover we say that $X$ is {\it left introverted} if, further
\begin{itemize}
\item[(iii)] $M\mult f\in X$ for $M\iin X^*$, with
$\norm{M\mult f}\leq\norm{M}\norm{f}$.
\end{itemize}
Notice that as an immediate consequence of (i), the evaluation functionals
$\eps_X(s)$, for $s\iin G$, are bounded; moreover
the family $\eps_X(G)=\{\eps_X(s)\}_{s\in G}$ is separating.
We observe that if we assume, in place of (iii),
only that $M\mult X\subset X$ for a fixed $M\iin X^*$
then it is automatic that the introversion operator $f\mapsto M\mult f$ is bounded.
Indeed, we appeal to the closed graph theorem:
if $\lim_{n\to\infty} f_n=f$ and $\lim_{n\to\infty}M\mult f_n=g$, then
for any $s\iin G$ we have $g(s)=\lim_{n\to\infty}M\mult f_n(s)=
\lim_{n\to\infty}M(f_n\mult s)=M(f\mult s)=M\mult f(s)$ and hence
$g=M\mult f$.  Similarly, if $X^*\mult f\subset X$, for a fixed $f\iin X$,
then $M\mapsto M\mult f$  is bounded.
However, we are aware of no means by which to prove
that the map $(M,f)\mapsto M\mult f$ is contractive.

If $X$ and $Y$ are both left homogeneous
Banach spaces in $\fC\fB(G)$, we say $X\subset Y$ boundedly (contractively)
if $X$ is a subspace of $Y$, and the inclusion map $X\hookrightarrow Y$
is bounded (contractive).

\begin{proposition}\label{prop:xtox}
Let $X$ be a left introverted  homogeneous subspace of $\fC\fB(G)$. Then:

{\bf (i)} $X^*$ is a right dual Banach algebra under the left Arens product;

{\bf (ii)} $\fX=\wbar{X}^{\unorm{\cdot}}$ is left introverted
with $X\subset\fX$ contractively; and

{\bf (iii)} if $Y$ is another left introverted homogeneous Banach space in
$\fC\fB(G)$, then $X\subset Y$ boundedly (contractively) if and only if
there is weak*-weak* continuous (contractive) operator $\Phi:Y^*\to X^*$ such that
$\Phi(\eps_Y(s))=\eps_X(s)$ for $s\iin G$.  The operator $\Phi$ is necessarily a homomorphism
(with respect to left Arens product).
If $X$ is a closed subspace of $Y$ then $\Phi$ is a quotient map.
\end{proposition}

\proof {\bf (i)}  This is standard, but short, so we include a full proof for completeness.
We have for $M,N\iin X^*$ that $\norm{(M\lap N)(f)}\leq
\norm{M}\norm{N\mult f}\leq \norm{M}\norm{N}\norm{f}$ for $f\iin X$, so
$\norm{M\lap N}\leq \norm{M}\norm{N}$.  Associativity remains:
if $L,M,N\in X^*$ and $f\in X$ then
\[
(M\lap N)\mult f(s)=M\lap N(f\mult s)=M(N\mult (f\mult s))=M((N\mult f)\mult s)
=M\mult(N\mult f)(s)
\]
for $s\iin G$ and hence
\[
L\lap(M\lap N)(f)=L((M\lap N)\mult f)=L(M\mult (N \mult f))=(L\lap M)(N\mult f)
=(L\lap M)\lap N(f).
\]
It is clear that $M\mapsto M\lap N$ is weak* continuous, so $X^*$ is a left dual Banach algebra.

{\bf (ii)}  If $m\in\fX^*$ and $f\in X$, then $m\mult f=M\mult f\in X$ where $M=m|_X$.
We note that $\unorm{m\mult f}\leq\norm{m}\sup_{s\in G}\unorm{f\mult s}
\leq\norm{m}\unorm{f}$.
By density of $X$ in $\fX$ and continuity of $f\mapsto m\mult f$ it follows that
$\fX$ is left introverted.

{\bf (iii)}  If $X\subset Y$ boundedly, we let $\Phi:Y^*\to X^*$ be the adjoint
of the inclusion map, which is the restriction map.  Then $\Phi$ intertwines $\eps_Y$ and $\eps_X$.
For $f\iin X$ and $N\iin Y^*$ we have
\[
N\mult f(s)=N(f\mult s)=\Phi(N)(f\mult s)=\Phi(N)\mult f(s)
\]
for $s\iin G$ so $N\mult f\in X$.  Then if, further,  $M\in X^*$, we obtain
\[
\Phi(M\lap N)(f)=M\lap N(f)=M(N\mult f)=\Phi(M)(\Phi(N)\mult f)=\Phi(M)\lap\Phi(N)(f).
\]
Conversely, if there exists a weak*-weak* continuous operator $\Phi:Y^*\to X^*$
which intertwines $\eps_X$ and $\eps_Y$, then the pre-adjoint $\vphi:X\to Y$
of $\Phi$ must satisfy $\vphi(f)(s)=\eps_Y(s)(\vphi(f))=\Phi(\eps_Y(s))(f)=\eps_X(s)(f)=f(s)$. Thus
$X\subset Y$ boundedly and $\vphi$ is the inclusion map.  If $X$ is a closed subspace
of $Y$, then $\Phi:Y^*\to X^*$, being the restriction map, is a quotient map
by the Hahn-Banach theorem.
\endpf

We consider a mild generalisation of Theorem \ref{theo:comparison} (i).

\begin{proposition}\label{prop:gx}
Let $X$ be a left introverted homogeneous subspace of $\fC\fB(G)$,
$\eps_X:G\to X^*$ be the evaluation map and
\begin{equation}\label{eq:sx}
G^X=\wbar{\eps_X(G)}^{w^*}\subset X^*.
\end{equation}
If $G^X$ is endowed with the weak* topology, then $(\eps_X,G^X)$ is
a right topological compactification of $G$.
\end{proposition}

We call $(\eps_X,G^X)$ the {\it $X$-compactification} of $G$.\medskip

\proof This proof is similar to an aspect of that of Theorem \ref{theo:comparison}, but short.
We first note that
\[
\eps_X(G)\subset Z_T(X^*)=\{M\in X^*:N\mapsto M\lap N\text{ is weak*-weak* continuous}\}.
\]
Indeed, if $s\in G$ and $N\in X^*$ we have for $f\iin X$ that
$\eps_X(s)\lap N(f)=\eps_X(s)(N\mult f)=N(f\mult s)$,
so $N\mapsto \eps_X(s)\lap N$ is weak*-weak* continuous.
Now if if $\chi,\chi'\in G^X$, we let
$\chi=\text{weak*-}\lim_\alp\eps_X(s_\alp)$ and
$\chi'=\text{weak*-}\lim_\beta\eps_X(t_\beta)$ for nets
$(s_\alp),(t_\beta)$ from $G$, and we have
\[
\chi\lap\chi'=\lim_\alp\eps_X(s_\alp)\lap\chi'=
\lim_\alp\lim_\beta\eps_X(s_\alp)\eps_X(t_\beta)=\lim_\alp\lim_\beta\eps_X(s_\alp t_\beta)
\in G^X.
\]
Being a subsemigroup of the right dual Banach algebra $X^*$ (see Proposition \ref{prop:xtox}
(i)), $G^X$ itself is right topological.
\endpf

We now consider two closed subalgebras of $\fC\fB(G)$ generated by
a left [right] introverted  homogeneous subspace $X$:
\[
\fA(X)=\wbar{\alg(X)}^{\unorm{\cdot}},\quad
\fE(X)=\wbar{\alg(X+\bar{X})}^{\unorm{\cdot}}\quad\aand\quad
\fE_1(X)=\fE(X)+\Cee 1
\]
where $\bar{X}$ denotes the space of complex conjugates of elements of $X$
and $\alg(X+\bar{X})$ the algebra generated by elements in $X$ and $\bar{X}$.
We note that $\fE_1(X)=\fE(X)$ if the latter is unital, and is the C*-unitization
otherwise.

\begin{theorem}\label{theo:cstarspec}
Let $X$ be a left introverted  homogeneous subspace of $\fC\fB(G)$.
Then

{\bf (i)} $\fE(X)$ and $\fE_1(X)$ are left m-introverted.

{\bf (ii)} $(\eps_{\fE_1(X)},\Phi_{\fE_1(X)})\cong(\eps_X,G^X)$ as
compactifications of $G$, and $\Phi_{\fE(X)}$ is homeomorphic
to $G^X\setdif\{0\}$; and

{\bf (iii)} $0\in G^X$ $\iff$ $1\not\in\fE(X)$.
\end{theorem}

\proof   If $\chi\in\Phi_{\fE(X)}$
let $\chi'=\chi|_X\in X^*$.  We have for a polynomial $p$ in $n+m$ variables with $p(0)=0$,
$f_1,\dots,f_n,g_1,\dots,g_m\iin X$, and $s\iin G$, that
\begin{align}\label{eq:charmultintop}
\chi\mult p(f_1,\dots,\bar{g}_m)(s)&=\chi(p(f_1\mult s,\dots,\wbar{g_m\mult s})) \notag \\
&=p(\chi'(f_1\mult s),\dots,\wbar{\chi'(g_m\mult s)})
=p(\chi'\mult f_1,\dots,\wbar{\chi'\mult g_m})(s)
\end{align}
so $\chi\cdot p(f_1,\dots,\bar{g}_m)\in\alg(X+\bar{X})$.
The introversion operator $f\mapsto \chi\mult f:\fE(X)\to\ell^\infty(G)$
is contractive, and takes a dense subspace of $\fE(X)$ into $\fE(X)$,
hence $\fE(X)$ is m-introverted.  If $\chi\in\Phi_{\fE_1(X)}$, then
$\chi(1)=1$ and the argument above, applied to an arbitrary polynomial,
shows that $\fE_1(X)$ is m-introverted.  Hence (i) is proved.

We shall prove (ii) and (ii) simultaneously.  Since $\fE_1(X)$ is a unital
C*-algebra of $\fC\fB(G)$, we have that $\Phi_{\fE_1(X)}=G^{\fE_1(X)}$.
The restriction map $\theta:\Phi_{\fE_1(X)}\to G^X$, $\theta(\chi)=\chi|_X$,
is weak*-weak* continuous which satisfies $\theta\comp\eps_{\fE_1(X)}(s)
=\eps_{\fE_1(X)}(s)|_X=\eps_X(s)$ for $s\iin G$, so $(\eps_X,G^X)\leq
(\eps_{\fE_1(X)},\Phi_{\fE_1(X)})$; in particular $\theta$ is surjective.
The map $\theta$ is injective
since each character is determined by its behavior on
$\alg(X+\bar{X})$ and hence on $X$.  Thus $\theta$ is a homeomorphism
and thus a compactification isomorphism.
Now if $1\in\fE(X)$ then $\Phi_{\fE(X)}=G^{\fE(X)}$.  Just as above,
the map $\theta:\Phi_{\fE(X)}\to G^X$ is surjective.  If
it were the case that $0\in G^X$, then for some $\chi\iin\Phi_{\fE(X)}$
$\chi|_X=0$, which would imply that $\chi(\alg(X+\bar{X}))=\{0\}$,
and imply that $\chi=0$, which is absurd.  If $1\not\in\fE(X)$,
then the unique character $\chi_\infty$ on $\fE_1(X)$ which
annihilates $\fE(X)$ satisfies $\theta(\chi_\infty)=0$, and thus, from the
surjectivity of $\theta$, $0\in G^X$.  In this case $\theta$
establishes a homeomorphism from  $\Phi_{\fE(X)}\cong\Phi_{\fE_1(X)}\setdif\{\chi_\infty\}$
onto $G^X\setdif\{0\}$. \endpf

We observe that it is possible that $G^X\setdif\{0\}$ is not a subsemigroup
of $G^X$.  If $G=\{o,e_1,e_2\}$ is the semilattice which is
generated by the relation $e_1e_2=o$, then
$\fX=\{f\in\fC(G):f(o)=0\}$ is an introverted subspace, in fact a subalgebra, for which
$G^\fX=\eps_\fX(G)\cong G$ and $0=\eps_\fX(o)\in G^\fX$.  Clearly
$G^\fX\setdif\{0\}\cong G\setdif\{o\}$.  A related example,
where $G$ is a group, is given in Example \ref{ex:sxnotsemigroup}, below.

The following is immediate from Theorem
\ref{theo:comparison}.

\begin{corollary}\label{cor:cstarspec2}
If $X,Y$ are two left homogeneous introverted subspaces of $\fC\fB(G)$
then $(\eps_X,G^X)\leq(\eps_Y,G^Y)$ if and only if $\fE(X)\subset\fE_1(Y)$.
In particular $(\eps_X,G^X)\cong(\eps_Y,G^Y)$ if and only if
 $\fE_1(X)=\fE_1(Y)$
\end{corollary}

We obtain an augmentation of Theorem \ref{theo:comparison} (i).
An open subset $U$ of a right topological semigroup $S$ has
{\it relatively proper right translations} if for every element $s\iin U$
and every compact subset $K$ of $U$, $Ks^{-1}\cap U$ is compact,
where $Ks^{-1}=\{t\in S:ts\in K\}$.

\begin{corollary}\label{cor:cstarspec1}
Let $\fX$ be a left translation invariant C*-subalgebra of $\fC\fB(G)$.
Then $\fX=\fC_0(G^\fX\setdif\{0\})\comp\eps_\fX$.
Moreover, $\fX$ is left  m-introverted if and only if $G^\fX$ is a right
topological semigroup for which $G^\fX\setdif\{0\}$ has relatively
proper right translations.
\end{corollary}

\proof We may and will assume that $1\not\in\fX$.
Let $\fX_1=\Cee1\oplus\fX$ be the
C*-unitisation of $\fX$.  From the theorem above we have
that $\Phi_{\fX_1}=G^{\fX_1}\cong G^\fX$
and $\Phi_\fX=G^\fX\setdif\{0\}$.
Since $\fX_1=\fC(G^\fX)\comp\eps_{\fX}$ and $0$, in $G^\fX$, corresponds to the
unique character which annihilates $\fX$, we have that
$\fX=\fC_0(G^\fX\setdif\{0\})\comp\eps_\fX$.

If $\fX$ is left m-introverted, then it is straighforward that $\fX_1$
is left m-introverted, so $G^\fX\cong G^{\fX_1}$ is a right topological
semigroup by Theorem \ref{theo:comparison} (i).  The condition
that $G^\fX\setdif\{0\}$ has relatively proper right translations is
equivalent to the condition that
\begin{equation}\label{eq:propertrans}
\chi\mult\fC_0(G^\fX\setdif\{0\})\subset\fC_0(G^\fX\setdif\{0\})
\text{ for every }\chi\iin G^\fX\setdif\{0\}.
\end{equation}
To see this, we require essentially the proof of \cite[Ex.\ 3.1.10]{berglundjm},
which we simply adapt to our situation.
If $U=G^\fX\setdif\{0\}$ has relatively proper right translations,
then for every $f\in\fC_0(U)$, $\eps>0$ and $\chi\in U$ we have
$\{\chi'\in U:\chi\mult f(\chi')\geq\eps\}=\{\chi'\in G^\fX:f(\chi')\geq\eps\}\chi^{-1}\cap U$
is compact, hence $\chi\mult f\in\fC_0(U)$.  Conversely, if
$U$ does not have relatively proper right translations, there is
compact $K\subset U$ and $\chi\in U$ for which the closed set $K\chi^{-1}\cap U$ is
not compact.  Then any compactly supported continuous $f:U\to[0,1]$
which satisfies $K\subset\{\chi':f(\chi')=1\}$ would satisfy that
$\chi\mult f\not\in\fC_0(U)$, hence (\ref{eq:propertrans}) cannot be satisfied.
Clearly, (\ref{eq:propertrans}) is necessary and sufficient for
$\fX=\fC_0(G^\fX\setdif\{0\})\comp\eps_\fX$ to be left m-introverted. \endpf

Let us consider some compactifications of $G$ which decompose with respect to
a second semigroup.  Suppose $H$ is locally compact right topological semigroup
for which there is a homomorphism $\eta:G\to H$ with dense range.
We say that $H$ has {\it proper} right translations if it has relatively proper
right translations on itself.
A left topological compactification $(\del,S)$ of $G$ is said to be an
{\it $(\eta,H)$-compactification} if there is a continuous $\theta:H\to S$ such that
$\theta\comp\eta=\del$; such a $\theta$ is necessarily a homomorphism.
Moreover, $(\del,S)$ is said to be a {\it regular $(\eta,H)$-compactification}
if $\theta$ is injective and open.
The following result is inspired by \cite[Lem.\ 4.1]{laul}
and \cite[Lem.\ 1]{ghahramanilaulosert}.

\begin{proposition}\label{prop:ghahlaulos}
Let $H$ be a locally compact noncompact right topological semigroup with proper right
translations, $\eta:G\to H$ be a homomorphism with dense range, and
$X$ be a left introverted homogeneous subspace of $\fC\fB(G)$.
Then $(\eps_X,G^X)$ is a regular $(\eta,H)$-compactification if and only if
$\fE(X)\supset\fC_0(H)\comp\eta$ and $\fC_0(H)\comp\eta$ is
an essential ideal in $\fE(X)$.

In this case there is a semigroup decomposition
\begin{equation}\label{eq:gxdecomp}
G^X=(G^X\setdif\theta(H))\sqcup\theta(H)
\end{equation}
where $\theta:H\to G^X$ is an injective, open homomorphism,
$\theta(H)$ is dense in $G^X$, and $G^X\setdif\theta(H)$ is a closed
ideal.
\end{proposition}

\proof From Theorem \ref{theo:cstarspec} (ii) we have that
$(\eps_X,G^X)\cong(\eps_{\fE(X)},G^{\fE(X)})$.  Hence it follows
Corollary \ref{cor:cstarspec1}, above, that
$\fE(X)=\fC_0(G^X\setdif\{0\})\comp\eps_X$.

If $(\eps_X,G^X)$ is a regular $(\eta,H)$-compactification
with open injective homomorphism $\theta$, then $\theta(H)\subset
G^X\setdif\{0\}$.  Indeed, if $\theta(H)\ni 0$, then $\{\theta^{-1}(0)\}$
would be an ideal in $H$, which would violate that $H$ has proper right
translations.  Since $\theta(H)$ is open and dense in $G^X\setdif\{0\}$,
it follows that $\fC_0(\theta(H))$ is an essential ideal in $\fC_0(G^X\setdif\{0\})$.
We thus have
\[
\fC_0(H)\comp\eta=\fC_0(H)\comp\theta^{-1}\comp\eps_X
=\fC_0(\theta(H))\comp\eps_X\subset\fC_0(G^X\setdif\{0\})\comp\eps_X=\fE(X)
\]
and $\fC_0(H)\comp\eta$ is an essential ideal in $\fE(X)$.  Conversely, if
$\fC_0(H)\comp\eta$ is an essential ideal in $\fE(X)=\fC_0(G^X\setdif\{0\})\comp\eps_X$,
then there is an open dense subset $U\subset G^X\setdif\{0\}$ such that
$\fC_0(H)\comp\eta=\fC_0(U)\comp\eps_X$.  Thus there is a homeomorphism
$\theta:H\to U$.  For $s\iin G$ we have that $\theta(\eta(s))=\eps_X(s)$,
i.e.\ for each $f\iin \fC_0(G^X\setdif\{0\})$ we have that
$f(\theta(\eta(s)))=f(\eps_X(s))$.
Hence $(\eps_X,G^X)$ is a regular $(\eta,H)$-compactification.

In the decomposition (\ref{eq:gxdecomp}), $\theta(H)$ is a dense open
subsemigroup by construction.  It remains to show that $G^X\setdif\theta(H)$
is an ideal.  We observe that $G^X\setdif\theta(H)=
\{\chi\in G^X:\chi|_{\fC_0(H)\comp\eta}=0\}$.  Hence
if $\chi\in G^X\setdif\theta(H)$ and $f\in\fC_0(H)$ then
for $s\iin G$, $\chi\mult (f\comp\eta)(s)=\chi(f\comp\eta(s))=0$,
so $\chi\mult(f\comp\eta)=0$.  Thus if $\chi'\in G^X$ and
$\chi\in G^X\setdif\theta(H)$ we have for $f\in\fC_0(H)$,
$\chi'\lap\chi(f)=\chi'(\chi\mult f)=0$, so $\chi'\lap\chi\in G^X\setdif\theta(H)$.
\endpf

The non-self-adjoint situation presents more complications than the self-adjoint one.
Since a non-self-adjoint uniform algebra $\fA$ on $G$ --- i.e.\
a uniformly closed subalgebra of $\fC\fB(G)$ ---
may admit spectrum larger than $G^\fA\setdif\{0\}$, the fact that the closure of the spectrum
is a semigroup must be checked.  Recall that our definition of left
m-introversion applies to all closed subalgebras of $\fC\fB(G)$.


\begin{proposition}\label{prop:phiasg}
Suppose $\fA$ is a closed left m-introverted subalgebra of $\fC\fB(G)$.
Then $\Phi_\fA\cup\{0\}$ is a semigroup under left  Arens product.
If $\fA$ is unital, then $\Phi_\fA$ itself is a semigroup.
\end{proposition}

As with non-unital self-adjoint algebras, we cannot expect that
$\Phi_\fA$ is a semigroup when $1\not\in\fA$. \medskip

\proof If $\chi\in\Phi_\fA$, $f,g\in \fA$ and $s\in G$ then similarly as in
(\ref{eq:charmultintop}) we have
\[
\chi\mult (fg)(s)=\chi(f\mult s\,g\mult s)=\chi(f\mult s)\chi(g\mult s)
=\chi\mult f\chi\mult g(s).
\]
Thus if we also have $\chi'\iin\Phi_\fA$ then
\[
\chi\lap\chi'(fg)=\chi(\chi'\mult(fg))=\chi(\chi'\mult f\,\chi'\mult g)
=\chi(\chi'\mult f)\chi(\chi'\mult g)=\chi\lap\chi'(f)\chi\lap\chi'(g)
\]
so $\chi\lap\chi'\in\Phi_\fA\cup\{0\}$.
If $\fA$ is unital, then for $\chi\in\Phi_\fA$ we have $\chi\mult 1=1$
and it follows for $\chi,\chi'\in\Phi_\fA$ that $\chi\lap\chi'(1)=1\not= 0$.
\endpf




If $\fA$ is a closed
subalgebra of a commutative C*-algebra, then a {\it boundary}
is any closed subset $B\subset\Phi_\fA$ such that $f\mapsto\hat{f}|_B:
\fA\to\fC_0(B)$ is an isometry ($f\mapsto\hat{f}$ is the Gel'fand transform).
The {\it \v{S}ilov boundary} is given by $\partial_\fA=\bigcap\{B:B\text{ is a boundary for }
\fA\}$, and is itself a boundary.  See the texts
\cite{kaniuthB,palmer} for details including the case that $\fA$ is non-unital.
In many cases where $\fA$ is a uniform algebra on $G$,
$\partial_\fA$ gives us a means of recovering $G^\fA$.

\begin{theorem}\label{theo:algspec}
Let $X$ be a left introverted homogeneous subspace of
$\fC\fB(G)$.

{\bf (i)} The algebra $\fA(X)$ is left m-introverted.
The map $\chi\mapsto\chi|_X:\Phi_{\fA(X)}\cup\{0\}\to X^*$
is a semigroup homomorphism and homeomorphism onto its range,
which takes $G^{\fA(X)}$ onto $G^X$.  Moreover $G^{\fA(X)}\setdif\{0\}$ is
closed in $\Phi_{\fA(X)}$.

{\bf (ii)}  We have $G^{\fA(X)}\setdif\{0\}\supset \partial_{\fA(X)}$.
Moreover, $\partial_{\fA(X)}$ is compact if $1\in\fE(X)$.

{\bf (iii)} If $G$ contains a dense subgroup $G_0$,
then $\wbar{\partial_{\fA(X)}}^{w*}$ is an ideal in $G^{\fA(X)}$.
Moreover, if $\eps_{X}(G_0)\cap \partial_{\fA(X)}\not=\varnothing$,
then $\partial_{\fA(X)}=G^{\fA(X)}\setdif\{0\}$.
\end{theorem}

The conclusion that $\wbar{\partial_{\fA(X)}}^{w*}$ is an ideal in $G^{\fA(X)}\cong G^X$
is false if we do not assume the existence of $G_0$.
See Example \ref{ex:algx} (iii), below, for this, and further illustrations.
We do not know if the condition
$\eps_{X}(G_0)\cap \partial_{\fA(X)}\not=\varnothing$
is automatic if the existence of $G_0$ is  assumed.
However, the latter condition is automatic if $G$ itself is
of a compact group; see, for example,
the proof of \cite[4.2.2]{grigoryant}.

\medskip
\proof That $\fA(X)$ is left m-introverted follows a calculation similar
to (\ref{eq:charmultintop}).  Hence it is immediate from Proposition \ref{prop:phiasg}
that $\Phi_{\fA(X)}\cup\{0\}$ is a subsemigroup of $\fA(X)^*$.  A simple modification
of the proof of Theorem \ref{theo:cstarspec} (ii), above, shows that
$\chi\mapsto\chi|_X:\Phi_{\fA(X)}\cup\{0\}\to\Phi_{\fA(X)}|_X\cup\{0\}$
is a homeomorphism which takes $G^{\fA(X)}$ onto $G^X$.  Since $G^{\fA(X)}$
is compact, it is closed in $\Phi_{\fA(X)}\cup\{0\}$, hence
$G^{\fA(X)}\setdif\{0\}$ is  closed in $\Phi_{\fA(X)}$.  The proof of
Proposition \ref{prop:xtox} (iii) shows that this restriction map is a semigroup
homomorphism.  Hence we have (i).

Since $\fA(X)$ generates the C*-algebra $\fE(X)$,
$G^{\fA(X)}\setdif\{0\}=(G^{\fE(X)}\setdif\{0\})|_{\fA(X)}=\Phi_{\fE(X)}|_{\fA(X)}$
is a boundary for $\fA(X)$, so $\partial_{\fA(X)}\subset G^{\fA(X)}\setdif\{0\}$.
Moreover, if $1\in\fE(X)=\fE(\fA(X))$, then by Theorem \ref{theo:cstarspec} (iii)
we have that $0\not\in G^{\fA(X)}$.
Hence $G^{\fA(X)}=G^{\fA(X)}\setdif\{0\}$ is compact, hence so too
must be the closed subset $\partial_{\fA(X)}$.  Hence we have (ii).

Now we consider (iii).
If $t\in G_0$ then $f\mapsto \eps_{\fA(X)}(t)\mult f$ and
$f\mapsto \eps_{\fA(X)}(t^{-1})\mult f$ are mutually inverse contractions, and hence
isometries.  We let for $t\iin G_0$ and $\chi\iin\Phi_{\fA(X)}$,
$t\mult\chi=\eps_{\fA(X)}(t)\lap\chi$.  Since $\eps_{\fA(X)}(G_0)\subset Z_T(G^{\fA(X)})$,
$\Phi_{\fA(X)}$ is a topological $G_0$-space.
If $B\subset\Phi_{\fA(X)}$ is any closed boundary, then $t\mult B$ must also be a
closed boundary for any $t\in G_0$.  Hence, by minimality, $t\cdot\partial_{\fA(X)}\supseteq
\partial_{\fA(X)}$, and it follows that $t\mult \partial_{\fA(X)}
=\partial_{\fA(X)}$.  Since $\eps_{\fA(X)}(G_0)$ is dense in $G^{\fA(X)}$
it follows that for $\chi\iin G^{\fA(X)}$, $\chi\lap \partial_{\fA(X)}\subset
\wbar{\partial_{\fA(X)}}^{w^*}$. Furthermore, if
$\eps_{X}(G_0)\cap \partial_{\fA(X)}\not=\varnothing$, then
$\eps_{X}(G_0)\subset \partial_{\fA(X)}$, and it follows that $\partial_{\fA(X)}$
is dense in $G^{\fA(X)}$.  Hence $\partial_{\fA(X)}=G^{\fA(X)}\setdif\{0\}$.
\endpf

We consider some properties associated with involutive semigroups.
Suppose now that $G$ has a continuous involution $s\mapsto s^*$,
i.e.\ $(s^*)^*=s$ and $(st)^*=t^*s^*$.
For $f\in\fC\fB(G)$ we define $f^*(s)=\wbar{f(s^*)}$.
A left [right] homogeneous subspace $X$ of $\fC\fB(G)$ will be called
{\it involutive} if it is closed under the involution $f\mapsto f^*$
and the involution is isometric on $X$.  In this case if
$M\in X^*$ (here $X^*$ still denotes the dual space), then we define
$M^*\in X^*$ by $M^*(f)=\wbar{M(f^*)}$.

\begin{proposition}\label{prop:involution}
Suppose $G$ admits a continuous involution and
$X$ is a left introverted involutive homogeneous subspace of $\fC\fB(G)$.

{\bf (i)}
The map  $M\mapsto M^*$ is a conjugate-linear, weak*-weak* continuous,
isometric involution on $X^*$.

{\bf (ii)}  The space $X$ is introverted.
On $X^*$ we have $(M\lap N)^*=N^*\rap M^*$ and $G^X$ is $*$-closed.

{\bf (iii)} The algebras
$\fE(X)$, $\fE_1(X)$, $\fA(X)$ and $\fA_1(X)$ are involutive
and $\Phi_{\fA(X)}$ is $*$-closed.
\end{proposition}

\proof The proof of (i) is straightforward.  We prove (ii).
Let $f\in X$ and $M\in X^*$.  Then for $s,t\iin S$ we have
$(s\mult f)^*(t)=\wbar{s\mult f(t^*)}=\wbar{f(t^*s)}=f^*(s^*t)=f^*\mult s^*(t)$ so
\[
f\mult M(s)=M(s\mult f)=\wbar{M^*(f^*\mult s^*)}=\wbar{M^*\mult f^*(s^*)}
=(M^*\mult f^*)^*(s)
\]
and hence $X$ is also right introverted, hence introverted.  Now it
follows for $M,N\iin X^*$ and $f\iin X$ that
\[
(M\lap N)^*(f)=\wbar{(M\lap N)(f^*)}=\wbar{M(N\mult f^*)}=M^*(f\mult N^*)
=N^*\rap M^*(f).
\]
Finally, if $\chi\in G^X$, then $\chi=\lim_\alp\eps_X(s_\alp)$ for some
net $(s_\alp)\subset G$, so we have for $f\iin X$
\[
\chi^*(f)=\wbar{\chi(f^*)}=\lim_\alp \wbar{f^*(s_\alp)}=\lim_\alp f(s_\alp^*)
\]
so $\chi^*=\text{weak*-}\lim_\alp \eps_X(s_\alp^*)\in G^X$ too.

We prove (iii).
If $p$ is a polynomial in $n$ variables with $p(0)=0$,
let $\bar{p}$ denote that same polynomial with
conjugated coefficients.  Now if $f_1,\dots,f_n\in X$ [or are
in $X\cup\bar{X}$], then $p(f_1,\dots,f_n)^*=\bar{p}(f_1^*,\dots,f_n^*)$ remains
in $\alg(X)$ [respectively $\alg(X+\bar{X})$].  Hence it follows from
continuity of the involution that $\fA(X)$ and $\fE(X)$ are involutive.
Clearly $1^*=1$ so $\fE_1(X)$ and $\fA_1(X)$ are involutive too.
If $\chi\in\Phi_{\fA(X)}$
then for $f,g\iin\fA(X)$, $\chi^*(fg)=\wbar{\chi(f^*g^*)}=\wbar{\chi(f^*)\chi(g^*)}
=\chi^*(f)\chi^*(g)$.  Clearly $\chi^*\not=0$ since $\chi\not=0$, so
$\chi^*\in\Phi_{\fA(X)}$ too.
\endpf

\subsection{On semi-topological compactifications}\label{ssec:semitopological}
We specialise our analysis above to semi-topological compactifications.
It is well known that $(\del,S)$ is a semitopological compactification of $G$
if and only if $(\del,S)$ is a factor of the weakly almost periodic compactification
$(\eps_{\fW\fA\fP},G^{\fW\fA\fP})$ associated to the weakly almost periodic functions $\fW\fA\fP(G)$;
or, equivalently, if and only if $\fC(S)\comp\del\subset\fW\fA\fP(G)$; see
\cite[\S 4.2]{berglundjm}, for example.

We summarise and build upon results due mainly to Glicksberg~\cite{glicksberg}
following Grothendieck~\cite{grothendieck} to prove a well-known result; see
\cite[4.2.7]{berglundjm}, for example.
We reprove this to demonstrate how
these properties amount to little more than properties of convolutions of measures.

\begin{theorem}\label{theo:glicksberg}
{\bf (i)}  Let $(\del,S)$ be a semitopological compactification of $G$.
Then $\fC(S)\comp\del$ is introverted and
on $(\fC(S)\comp\del)^*\cong\meas(S)$ the left and right Arens
products coincide and are given by convolution:
\begin{equation}\label{eq:conv}
\mu\con\nu(f)=\iint_{S\times S} f(\chi'\chi)d\nu(\chi)d\mu(\chi')
=\iint_{S \times S}f(\chi'\chi)d\mu(\chi')d\nu(\chi)
\end{equation}
for $\mu,\nu\in\meas(S)$ and $f\iin\fC(S)$

{\bf (ii)} Let $\fX$ be a closed, translation invariant subspace
of $\fW\fA\fP(G)$.  Then $\fX$ is introverted and Arens regular.
\end{theorem}

\proof Let $\mu,\nu\in\meas(S)\cong\fC(S)^*$ and $f\in\fC(S)$.  Then
for $s\iin G$ we have
\[
\nu\mult(f\comp\del)(s)=\nu((f\comp\del)\mult s)
=\int_S f(\del(s)\chi)d\nu(\chi)
\]
and, similarly, $(f\comp\del)\mult\mu(s)=\int_S f(\chi'\del(s))d\mu(\chi')$.
Thanks to \cite[1.2]{glicksberg} we have that $\nu\mult f$ and $f\mult\mu$
given by
\[
\nu\mult f(\chi)= \int_S f(\chi\chi')d\nu(\chi')\aand
f\mult\mu(\chi)= \int_S f(\chi'\chi)d\mu(\chi')
\]
are elements of $\fC(S)$, hence $\fC(S)\comp\del$ is introverted.
Thus, the Fubini theorem \cite[3.1]{glicksberg} shows that
\begin{align*}
\mu\lap\nu(f\comp\del)&=\mu(\nu\mult (f\comp\del))
=\iint_{S\times S} f(\chi'\chi)d\nu(\chi)d\mu(\chi')  \\
&=\iint_{S \times S}f(\chi'\chi)d\mu(\chi')d\nu(\chi)=\nu((f\comp\del)\mult\mu)
=\mu\rap\nu(f\comp\del)
\end{align*}
Hence left and right Arens products coincide on measures as functionals
on $\fC(S)$.  Thus (i) is established.

We prove (ii).  By Theorem \ref{theo:cstarspec} (ii), $G^\fX\cong G^{\fE_1(\fX)}$.
Since $\fE_1(\fX)\subset\fW\fA\fP(G)$, it follows Theorem \ref{theo:comparison}
that $S=G^\fX$ is a semitopological semigroup.
We have that $\fX=\fF\comp\eps_{\fX}$ for some closed $\del(G)$-translation invariant
subspace $\fF$ of $\fC(S)$.  We first note that $\fF$ is also $S$-translation
invariant.  Given $s\iin S$ let $(t_\alp)$ be a net from $G$ so
$s=\lim_\alp\del(t_\alp)$.  Hence for $s'\iin S$ we have
$s\mult f(s')=f(s's)=\lim_\alp f(s'\del(t_\alp))
=\lim_\alp \del(t_\alp)\mult f(s)$, so
$s\mult f=\text{pointwise-}\lim_\alp \del(t_\alp)\mult f$.
By \cite[Theo.\ 5]{grothendieck}, also see \cite[A.2]{berglundjm},
$s\mult f=\text{weak-}\lim_\alp \del(t_\alp)\mult f$.  Hence it follows from the Hahn-Banach
theorem that $s\mult f\in\fF$.  A symmetric argument gives that
$f\mult s\in\fF$ too.

Now, let $f\in\fF$, $m\in\fF^*$ and let $\mu\iin\meas(S)$ be so $\mu|_\fF=m$.
We note that $m\mult f=\mu\mult f$.
We let $\mu=\text{weak*-}\lim_\alp\mu_\alp$, where each
$\mu_\alp$ is a finite linear combination of point masses of norm
not exceeding $\norm{\mu}_1$, which may be realised with aid of the
Krein-Milman theorem.
Since $\fF$ is translation invariant,
$\mu_\alp\mult f\in\fX$ for each $\alp$, and hence for $s\in S$
we have $\mu\mult f(s)=\mu(f\mult s)=\lim_\alp \mu_\alp(f\mult s)
=\lim_\alp \mu_\alp\mult f(s)$.
Thus $\mu\mult f=\text{pointwise-}\lim_\alp \mu_\alp\mult f$,  and, as above,
we deduce that $\mu\mult f\in\fF$.  By a symmetric argument we have that
$f\mult \mu\in\fF$ too.  It follows that $\fX=\fF\comp\eps_{\fX }$ is introverted.
Finally, since $\mu\mapsto\mu|_\fF$ is a quotient homomorphism
by Proposition \ref{prop:xtox} (iii),
Arens regularity of $\fC(S)$ passes to that of $\fF$, and hence to $\fX$.
\endpf

\begin{corollary}\label{cor:semitopalg}
If $X$ is a homogeneous subspace of $\fW\fA\fP(G)$ then

{\bf (i)}  $G^X$ is a semitopological semigroup;

{\bf (ii)} $\fE(X)$ and $\fA(X)$ are introverted;

{\bf (iii)} the Arens product on $\fE(X)^*\cong\meas(G^X\setdif\{0\})$
is given by convolution product: if $0\in G^X$ we let for
$f\iin\fC_0(G^X\setdif\{0\})$ and $\mu,\nu\iin \meas(G^X\setdif\{0\})$
\begin{equation}\label{eq:conv1}
\mu\con\nu(f\comp\eps_X)=\int_{G^X\setdif\{0\}}f(\chi)d(\mu\con\nu)(\chi)
=\iint_{G^X\times G^X}\til{f}(\chi\chi')d\til{\mu}(\chi)d\til{\nu}(\chi')
\end{equation}
where $\til{f}$ is the continuous extension of $f$ to $G^X$ satisfying
$\til{f}(0)=0$, and $\til{\mu},\til{\nu}\iin\meas(G^X)$ are any measures
which restrict on Borel subsets of $G^X\setdif\{0\}$ to $\mu$ and $\nu$; and

{\bf (iv)} $\Phi_{\fA(X)}\cup\{0\}$ is a semitoplogical semigroup under
convolution product.
\end{corollary}

\proof It follows from the proof of Theorem \ref{theo:cstarspec} that
$\fE_1(X)$ is a translation invariant unital C*-subalgebra of $\fW\fA\fP(G)$,
and hence $G^{\fE_1(X)}\cong G^X$ is a semi-topological semigroup.  Thus we obtain
(i). Part (ii) is immediate from (ii) of the theorem above.

By Theorem \ref{theo:cstarspec}, the convolution formula (\ref{eq:conv1})
of part (iii) extends (\ref{eq:conv}) only when $1\not\in\fE(X)$.  In this case
$f\mapsto\til{f}:\fC_0(G^X\setdif\{0\})\to\fC(G^X)$ is the canonical embedding
and we essentially appeal to Proposition \ref{prop:xtox} (iii).
We remark that one can select $\til{\mu}$ and $\til{\nu}$, in (\ref{eq:conv1})
to satisfy $\til{\mu}(\{0\})=0=\til{\nu}(\{0\})$.

Finally (iv) is immediate from Proposition \ref{prop:phiasg} and the theorem above.
\endpf

It may fail that $X$, above is itself introverted.

\begin{example}
{\bf (i)} Let $G=\mathrm{F}_\infty$, the free group on countably many generators, and
consider the space $X=\ccoef^{cb}(G)$, the norm-closure of
$\falg$ in the completely bounded multipliers $\cbmg$.
See \cite{decanneireh} for the definition of $\cbmg$
and its isometric predual $\queg$,
and \cite{forrestrs} for more on $\ccoef^{cb}(G)$.
It is noted in \cite{xu} that $\cbmg\subset \fW\fA\fP(G)$.
If $\ccoef^{cb}(G)$ were itself introverted
then by Proposition \ref{prop:xtox},
its dual space $\ccoef^{cb}(G)^*$ would be a Banach subalgebra contained in
$\vn(G)\cong\falg^*$ containing the algebra $\lam_1(\ell^1(G))$.
The weak amenability property of $G$ (see \cite{decanneireh}) tells us that
$\ccoef^{cb}(G)$ is weak*-dense in $\cbmg$, hence the adjoint of
the inclusion map $\ccoef^{cb}(G)\hookrightarrow\cbmg$ takes
$\queg$ isometrically onto the closed subspace
generated $\lam_1(\ell^1(G))$ in $\ccoef^{cb}(G)^*$.
However computations of Haagerup \cite{haagerup} (see \cite[Remark 3.2]{pisier})
show that $\queg$ is not a Banach algebra with respect to this product.
Thus it cannot be the case that $\ccoef^{cb}(G)$ is introverted.

{\bf (ii)}  Let $G$ be any infinite discrete group.  Then $\ell^1(G)$
is a homogeneous subspace of $\fC_0(G)\subset\fW\fA\fP(G)$, but is
not introverted.  Indeed the constant function
$1\iin\ell^\infty(G)$ satisfies $1\mult f(s)=\sum_{t\in G}f(t)$ for each
$s$, i.e.\ is constant, and is generally not an element of $\ell^1(G)$.
\end{example}

If $G$ is an involutive semigroup, a
right [left] topological compactification
$(\del,S)$ of $G$ is called {\it involutive} if there is a continuous
involution $x\mapsto x^*$ on $S$ such that $\del(s^*)=\del(s)^*$
for $s\iin G$.

\begin{proposition}\label{prop:semitopinvolution}
Suppose $G$ admits a continuous involution.

{\bf (i)} If a right topological semigroup $S$ admits a continuous
involution, then it is semitopological.

{\bf (ii)} If $(\del,S)$ is an involutive right topological compactification, then
$S$ is semitopological, and $\fC(S)\comp\del$ is involutive in the sense
of Proposition \ref{prop:involution}.

{\bf (iii)} The weakly almost periodic compactification
$(\eps_{\fW\fA\fP},G^{\fW\fA\fP})$ is an involutive compactification.

In particular, $(\eps_{\fW\fA\fP},G^{\fW\fA\fP})$
is the universal involutive compactification of $G$.
\end{proposition}

\proof  We have for $s\iin S$ that $t\mapsto t^*s^*$ is continuous on
$S$, and hence $t\mapsto st=(t^*s^*)^*$ is continuous on $S$, so
$s\in Z_T(S)$.  Thus we
have (i).  Part (ii) is immediate from (i), and the fact that for $f\iin\fC(S)$,
$(f\comp\del)^*=f^*\comp\del$, where $f^*$ is clearly in $\fC(S)$.

To see (iii) we need only show that $\fW\fA\fP(G)$ is $*$-closed,
from which point we may appeal to Proposition \ref{prop:involution} (i) and (ii).
It follows from \cite[4.2.3]{berglundjm},
$f\in\fW\fA\fP(G)$ if and only if $G\mult f$, or equivalently $f\mult G$,
is relatively weakly compact in $\fC\fB(G)$.
Since $(s\mult f^*)^*=f\mult s^*$
for $s\iin G$, it follows that $G\mult f^*$ is relatively weakly compact
exactly when $f\mult G$ is.

Since $(\eps_{\fW\fA\fP},G^{\fW\fA\fP})$ is universal amongst
all semi-topological compactifications, it follows from (i) that it dominates
any involutive compactification. Hence $(\eps_{\fW\fA\fP},G^{\fW\fA\fP})$
is the universal involutive compactification.
\endpf




The following should be compared to results in \cite{farhadig}.

\begin{corollary}\label{cor:gxinvolutive}
Suppose $G$ admits a continuous involution and $X$ is a left introverted
homogeneous subspace of $\fC\fB(G)$.  Then $G^X$ admits an involution
by which $(\eps_X,G^X)$ is
an involutive compactification if and only if $X\subset\fW\fA\fP(G)$
and $\fE(X)$ is involutive in the sense of Proposition \ref{prop:involution}.
\end{corollary}

Curiously, we need not assume that $X$ itself is involutive in the sense of
Proposition \ref{prop:involution}, since properties of $G^X$ are
determined by the structure of $\fE(X)$.  Moreover, it does not appear
to be the case that having $\fE(X)$ involutive necessarily implies
that $X$ itself must be, though we have no examples to suggest otherwise.
\medskip

\proof We have $\fE_1(X)=\fC(G^X)\comp\eps_X$ by
Theorem \ref{theo:comparison} (i) and Theorem \ref{theo:cstarspec} (ii).
It then follows Proposition \ref{prop:semitopinvolution}
and Theorem \ref{theo:comparison} (ii) that if $(\eps_X,G^X)$
is involutive then $\fE_1(X)$, and hence $X$, is contained in
$\fW\fA\fP(G)$.  Also $\fE_1(X)$ is clearly involutive; since
$1^*=1$, necessarily $\fE(X)$ is involutive too.

Conversely, if $X\subset\fW\fA\fP(G)$ and $\fE(X)$ is involutive,
then $\fE_1(X)$ is involutive and contained in $\fW\fA\fP(G)$.
Hence by Proposition \ref{prop:involution} (ii),
and Theorem \ref{theo:comparison} (ii) $(\eps_X,G^X)$
is involutive as well.
\endpf

\begin{example}\label{ex:ucbandluc}
We note that if $G$ is a group,
the space of uniformly continuous bounded functions $\fU\fC\fB(G)$
is $*$-closed where $s^*=s^{-1}$ for $s\iin G$.
In general $G^{\fU\fC\fB}=G^{\fU\fC\fB(G)}$ is not an involutive semigroup
in our sense, i.e.\ the involution is not continuous.
We note that $G$ has the small invariant neighbourhood [SIN] property if and only
if $\fU\fC\fB(G)=\fL\fU\fC(G)$; see \cite[(4.14)(g)]{hewittrI}, for example.
Thus if $G$ is a locally compact non-compact [SIN]-group then $\fU\fB\fC(G)$
is not Arens regular, since in this case
neither is $\fL\fU\fC(G)$; see \cite{lau} or \cite{neufang} for example.
\end{example}

\subsection{Compactifications which are semigroups of contractions on Hilbert spaces}
\label{ssec:CH}
If $\fH$ is a Hilbert space, let $\fB(\fH)_{\norm{\cdot}\leq 1}$ denote the
weak* (weak operator) compact semitopological semigroup of
linear contractions on $\fH$.  A semigroup of Hilbertian contractions
is any subsemigroup $S\subset\fB(\fH)_{\norm{\cdot}\leq 1}$.
We say
\begin{center}
(CH)\qquad
\parbox{3.7in}{
a compactification $(\del,S)$ of $G$ has property (CH) if $S$
is isomorphic to a weak*-closed semigroup of Hilbertian contractions.}
\end{center}
We see that there is a universal (CH)-compactification.

\begin{theorem}\label{theo:chuniversal}
There is a (CH)-compactification $(\eps_{\fC\fH},G^{\fC\fH})$ which is universal
in the sense that every (CH)-compactification of $G$ is a factor
of $(\eps_{\fC\fH},G^{\fC\fH})$.
\end{theorem}

\proof Thanks to \cite[3.3.4]{berglundjm} it suffices to verify that the class
of (CH) compactifications is closed under subdirect products.
If $\{(\del_i,S_i)\}_{i\in I}$ is a set of (CH) compactifications, where each $S_i$ is
isomorphic to a weak*-closed subsemigroup of $\fB(\fH_i)_{\norm{\cdot}\leq 1}$, then
$P=\prod_{i\in I}S_i$ is isomorphic to a weak*-closed subsemigroup of
$\prod_{i\in I}\fB(\fH_i)_{\norm{\cdot}\leq 1}\til{\subset}
\fB(\ell^2\text{-}\bigoplus_{i\in I}\fH_i)_{\norm{\cdot}\leq 1}$.
Thus if $\del:G\to P$ is given by $\del(s)=(\del_i(s))_{i\in I}$, and
$S=\wbar{\del(G)}^{w^*}$,
then $(\del,S)$ is the subdirect product $\bigvee_{i\in I}(\del_i,S_i)$ of this system,
in the sense of \cite[3.2.5]{berglundjm}.  Clearly $(\del,S)$ is
a (CH) compactification.  \endpf

As a consequence of Theorem \ref{theo:comparison}, we see that
if $\fC\fH(G)=\fC(G^{\fC\fH})\comp\eps_{\fC\fH}$, then $\fC\fH(G)\supset\fC(S)\comp\del$
for every (CH)-compactification $(\del,S)$ of $G$.  What is not clear
is whether or not the class of (CH)-compactifications is closed under
homomorphism.  Thus we call a compactification $(\del,S)$ of $G$ an
(FCH)-compactification if it is a factor of a (CH)-compactification, hence a factor
of $(\eps_{\fC\fH},G^{\fC\fH})$.  It is immediate from Theorem \ref{theo:comparison} (ii)
that  $(\del,S)$ is an (FCH)-compactification if and only if $\fC(S)\comp\del
\subset\fC\fH(G)$.

Let us see that each concrete (CH) compactification is associated with a particular
operator algebra and an associated homogeneous Banach space in $\fC\fB(G)$.
This is modeled closely on \cite[(2.2)]{arsac}.

\begin{theorem}\label{theo:arsac}
Let $\del:G\to(\fB(\fH)_{\norm{\cdot}\leq 1},w^*)$ be a continuous homomorphism
and $\opalg_\del=\wbar{\spn\del(G)}^{w^*}$.  Then there is an introverted homogeneous
Banach space $\ccoef_\del$ in $\fC\fB(G)$ such that $\ccoef_\del^*\cong
\opalg_\del$.  Moreover $\opalg_\del$ is a subalgebra of $\fB(\fH)$,
and the Arens products on $\opalg_\del$ coincide and are exactly the
operator products.
\end{theorem}

\proof The unique predual of $\fB(\fH)$ is given by the projective tensor
product $\fH\what{\otimes}\bar{\fH}$, via the identification
$T(\xi\otimes\bar{\eta})=\inprod{T\xi}{\eta}$.
Define $E_\del:\fH\what{\otimes}\bar{\fH}\to\fC\fB(G)$ by
\[
E_\del x(s)=\sum_{j=1}^\infty\inprod{\del(s)\xi_j}{\eta_j}
\]
for $x=\sum_{j=1}^\infty\xi_j\otimes\bar{\eta}_j\iin\fH\what{\otimes}\bar{\fH}$
and $s\iin G$.  If we let
$\ccoef_\del=\mathrm{ran}E_\del$ be given the norm by which
$E_\del$ is a quotient map, then $\ccoef_\del$ is a Banach space
and we have for $\eps>0$ and $x$ as above with
$\sum_{j=1}^\infty\norm{\xi_j}\norm{\eta_j}<\norm{x}+\eps<\norm{E_\del x}+2\eps$
that for $s\iin G$
\[
|E_\del x(s)|\leq\sum_{j=1}^\infty|\inprod{\del(s)\xi_j}{\eta_j}|
\leq\norm{E_\del x}+2\eps
\]
so $\unorm{E_\del x}\leq\norm{E_\del x}$, and $\ccoef_\del$ is
indeed a subspace of $\fC\fB(G)$.
The bipolar theorem shows that $\opalg_\del=(\ker E_\del)^\perp$ and hence
it follows from the Hahn-Banach theorem that
$\ccoef_\del^*\cong(\ker E_\del)^\perp=\opalg_\del$.  The duality
relation may be realised by $N(\inprod{\del(\cdot)\xi}{\eta})
=\inprod{N\xi}{\eta}$ for $N\iin\opalg_\del$ and $\xi,\eta\iin\fH$.

We verify that $\ccoef_\del$ is introverted. First note that for  $\xi,\eta\iin\fH$
we have
\[
\inprod{\del(\cdot)\xi_j}{\eta_j}\mult s=\inprod{\del(\cdot)\del(s)\xi_j}{\eta_j}
\aand s\mult \inprod{\del(\cdot)\xi_j}{\eta_j}=\inprod{\del(\cdot)\xi_j}{\del(s)^*\eta_j}
\]
from which it follows that $\ccoef_\del$ is translation invariant.
Now, for $N\iin\opalg_\del$ and
$\xi,\eta\in\fH$ we have that
\[
N\mult\inprod{\del(\cdot)\xi}{\eta}(s)=N(\inprod{\del(\cdot)\xi}{\del(s)^*\eta})
=\inprod{\del(s)N\xi}{\eta}.
\]
Hence, if $N\iin\opalg_\del$, and $x\iin\fH\what{\otimes}\bar{\fH}$ as above,
we have
\[
\norm{N\mult(E_\del x)}=
\norm{\sum_{j=1}^\infty\inprod{\del(\cdot)N\xi_j}{\eta_j}}
\leq \sum_{j=1}^\infty\norm{N\xi_j}\norm{\eta_j}
\leq\norm{N}(\norm{E_\del x}+2\eps).
\]
A similar calculation establishes that $\norm{(E_\del x)\mult N}\leq
\norm{N}(\norm{E_\del x}+2\eps)$, hence it follows that $\ccoef_\del$
is introverted.  In particular, setting $N=\del(s)$ we see that
$\ccoef_\del$ is homogeneous as well.

It is immediate that, $\opalg_\del$, being a weak*-closure
of a subalgebra of the dual Banach algebra $\fB(\fH)$, is itself closed under
operator multiplication.  Let us check the left Arens product.
If $M,N\in \opalg_\del$ and $\xi,\eta\iin\fH$ we have
\[
M\lap N(\inprod{\del(\cdot)\xi}{\eta}
=M(N\mult\inprod{\pi(\cdot)\xi}{\eta})=M(\inprod{\del(\cdot)N\xi}{\eta})
=\inprod{MN\xi}{\eta}.
\]
The computations for the right Arens product are similar.  \endpf

We note that if $G$ is an involutive semigroup and
$\del:G\to(\fB(\fH)_{\norm{\cdot}\leq 1},w^*)$ is
a continuous $*$-homomorphism, then $\opalg_\del$ is a $*$-subalgebra
of $\fB(\fH)$, in fact a von Neumann algebra if $\del$ is non-degenerate.
Also, $\ccoef_\del$ is involutive in the sense of Section \ref{ssec:homogeneousspaces},
by a calculation similar to that for (\ref{eq:inversioninvariance}).  We give
a mild refinement of Proposition \ref{prop:xtox}.  We maintain the
convention of Section \ref{ssec:homogeneousspaces} of stating results
only for left introverted spaces and right topological semigroups.

\begin{proposition}\label{prop:xtovn}
Let  $\del:G\to(\fB(\fH)_{\norm{\cdot}\leq 1},w^*)$ be
a continuous homomorphism,
and $X$ a left introverted homogeneous Banach space of
$\fC\fB(G)$.  Then
$\ccoef_\del\subset X$ boundedly (contractively)
if and only if there is a weak*-weak* continuous (contractive) homomorphism
$\del_X:X^*\to\opalg_\del$ for which $\del_X(\eps_X(s))=\del(s)$ for $s\iin G$.
Moreover, if $G$ is an involutive semigroup, $\del$ is a $*$-homomorphism, and
$X$ is involutive, then $\del_X(M^*)=\del_X(M)^*$ for
$M\iin X^*$.
\end{proposition}

\proof The first statement is immediate from Proposition \ref{prop:xtox} (iii)
and Theorem \ref{theo:arsac}.


Further if $X$ is involutive and $\del$ is a $*$-homomorphism, then for
$M\in X^*$ we have for $\xi,\eta\iin\fH$, using (\ref{eq:inversioninvariance}), that
\[
\inprod{\del_X(M^*)\xi}{\eta}=\wbar{M(\inprod{\del(\cdot)\eta}{\xi})}
=\inprod{\del_X(M)^*\xi}{\eta}.
\]
Notice that $\del_X$ is a homomorphism with respect to either Arens product
and $\del_X(M\lap N)=\del_X(M\rap N)$.  This is to be expected as
elements of $\ccoef_\del$ are weakly almost periodic functions.
\endpf

We shall see in Corollary \ref{cor:universal2} that there
are compact semitopological involutive semigroups
$G$ for which $G^{\fC\fH}\not= G$.
Hence for such semigroups no analogue of the Gelfand-Raikov theorem
(as stated in \cite[(22.12)]{hewittrI}, for example) holds.


\section{Eberlein compactifications}\label{sec:eberlein}

From this point forward, let $G$ denote a locally compact group.

\subsection{Basic theory}
If $\pi\in\Sig_G$ we let $G^\pi=\wbar{\pi(G)}^{w^*}\subset\vn_\pi
\cong\ccoef_\pi^*$.  The pair
\[
(\pi,G^\pi)
\]
is called the {\it $\pi$-Eberlein compactification} of $G$.  Moreover we set
\[
(\eps_\fE,G^\fE)=(\ome_G,G^{\ome_G})
\]
and simply call this the {\it Eberlein compactification} of $G$.
The representations $\tau_\pi$ and $\rho_\pi$ and spaces $\ccoef_\pi,\fE_\pi,
\ccoef(\pi),\fA(\pi),\ecoef(\pi)$ and $\fE(\pi)$ are defined in
Section \ref{ssec:basicspaces}.
We remark that in the notation of Section \ref{ssec:homogeneousspaces} we have
\[
\fA(\pi)=\fA(X),\;\fE(\pi)=\fE(X)\quad\text{and}\quad
\fE_1(\pi)=\fE_1(X)
\]
where $X=\ccoef_\pi$; in particular $(\pi,G^\pi)=(\eps_X,G^X)$.

\begin{proposition}\label{prop:gpiisspec}
If $\pi\in\Sig_G$, then $(\pi,G^\pi)$ is an involutive
semitopological compactification of $G$.
We have an equivalence of compactifications
\[
(\pi,G^\pi)\cong(\tau_\pi,G^{\tau_\pi})\cong(\rho_\pi,G^{\rho_\pi}).
\]
Moreover, $G^\pi$ is the Gelfand spectrum of the $\pi$-Eberlein algebra
$\fE_1(\pi)$, and $G^\pi\setdif\{0\}$ is that of $\fE(\pi)$.
\end{proposition}

\proof That $(\pi,G^\pi)$ is an involutive semitopological compactification can be verified
directly as involution of $\vn_\pi$ is weak*-weak* continuous.  However this
can also be deduced from results above.
Indeed, on $X=\ccoef_\pi$ we note that $f^*=\bar{\check{f}}$,
so $X$ is involutive by (\ref{eq:inversioninvariance}).  Thus
$X^*$ admits a linear involution by Proposition \ref{prop:involution}.
Finally, $X\subset\fW\fA\fP(G)$, by
\cite[Theo.\ 3.11]{burckel} or \cite[Theo.\ 11.3]{eberlein},
for example, so the involution on $G^X$ is a continuous
semigroup involution by Corollary  \ref{cor:gxinvolutive}.

We have that $(\pi,G^\pi)\cong(\eps_{\fE_1(\pi)},\Phi_{\fE_1(\pi)})$
by Theorem \ref{theo:cstarspec}; and from there we also identify the
spectrum of $\fE(\pi)$.  Since $\fE_1(\pi)=\fE_1(X)$ for
either of $X=\ccoef(\pi),\ecoef(\pi)$, the equivalence of each
of $(\tau_\pi,G^{\tau_\pi})$ and $(\rho_\pi,G^{\rho_\pi})$ with
$(\pi,G^\pi)\cong(\eps_X,G^X)$ holds by Corollary \ref{cor:cstarspec2}.
\endpf

Let us review some established methods of comparing representations
$\pi$, $\sig$ in $\Sig_G$.  We say
\begin{itemize}
\item[$\bullet$] $\pi$ is {\it quasi-contained} in $\sig$, $\pi\qc\sig$, if
$F_\pi\subset\ccoef_\sig$; and
\item[$\bullet$] $\pi$ is {\it weakly contained} in $\sig$, $\pi\wc\sig$, if
$F_\pi\subset\wcoef_\sig$.
\end{itemize}
We add two more useful comparisons:
\begin{itemize}
\item[$\bullet$] $\pi$ is {\it m-quasi-contained} in $\sig$, $\pi\mq\sig$, if
$F_\pi\subset\ccoef(\sig)$; and
\item[$\bullet$] $\pi$ is {\it \={m}-quasi-contained} in $\sig$, $\pi\cmq\sig$, if
$F_\pi\subset\ecoef(\sig)$.
\end{itemize}
Using our spaces we define notions of ``Eberlein containment'':
\begin{itemize}
\item[$\bullet$] $\pi$ is {\it strongly Eberlein contained} in $\sig$, $\pi\se\sig$, if
$F_\pi\subset\fE_\sig$;
\item[$\bullet$] $\pi$ is {\it m-Eberlein contained} in $\sig$, $\pi\me\sig$, if
$F_\pi\subset\fA(\sig)$;
\item[$\bullet$] $\pi$ is {\it  Eberlein contained} in $\sig$, $\pi\ec\sig$, if
$F_\pi\subset\fE(\sig)$; and
\item[$\bullet$] $\pi$ is {\it  weakly Eberlein contained} in $\sig$, $\pi\we\sig$, if
$F_\pi\subset\fE\fB(\sig)$,
\end{itemize}
where $\fE\fB(\pi)$ is defined in Section \ref{ssec:basicspaces}.
Each of the relations above is transitive with the exception of weak Eberlein containment;
see Example \ref{ex:rachman} (ii), below.
For each comparison $\pi\leq_\bullet\sig$, above, we obtain an equivalence
$\pi\cong_\bullet\sig$ in the case that $\sig\leq_\bullet\pi$ holds, as well.

The following gives the motivation for including
conjugation and multiplication properties in our definition of
Eberlein containment.  It is immediate from Corollary \ref{cor:cstarspec2} and
Proposition \ref{prop:gpiisspec}.

\begin{proposition}\label{prop:eberleincont}
If $\pi,\sig\in\Sig_G$ we have

{\bf (i)} $\pi\ec\sig\oplus 1$ if and only if $(\pi,G^\pi)\leq(\sig,G^\sig)$, and

{\bf (ii)} $\pi\oplus 1\cong_e\sig\oplus 1$ if and only if  $(\pi,G^\pi)\cong(\sig,G^\sig)$.
\end{proposition}

For $\pi$, $\sig$ in $\Sig_G$ our definitions provide
the following implications.
\begin{equation} \label{eq:diagram}
\xymatrix{
\pi\qc\sig \ar@{=>}[rr]\ar@{=>}[dd]\ar@{=>}[dr]
& & \pi\mq\sig \ar@{=>}[r]\ar@{=>}[d] & \pi\cmq\sig \ar@{=>}_{(1)}[d]  \\
& \pi\se\sig \ar@{=>}[r]& \pi\me\sig \ar@{=>}[r] &  \pi\ec\sig \ar@{=>}[d]  \\ 
\pi\wc\sig\ar@{=>}^{(2)}[rrr]  & & & \pi\we\sig
}
\end{equation}
We shall note, after Theorem \ref{theo:algcontone}, that this diagram simplifies 
significantly if $\pi$ is finite dimensional.
None of the implications above are equivalences, in general, showing that
these notions of containment are distinct.  We will show this, for implications 
(1) and (2), by way of some examples.

\begin{example} \label{ex:rachman}
For any non-compact $G$, $\foalg=\fsalg\cap\fC_0(G)$
is a closed, translation invariant, conjugate invariant subalgebra of $\fsalg$,
called the Rajchman algebra.  Hence, by \cite[(3.17)]{arsac} there is a representation
$\rho_0$ for which $\foalg=\ccoef_{\rho_0}=\ecoef(\rho_0)$.  In general
\begin{equation}\label{eq:rachmancont}
\lam\qc\rho_0\text{, in particular }\lam\cmq\rho_0.
\end{equation}
The well known fact that $\fE_\lam=\fC_0(G)$ gives
$\lam\cong_e\rho_0$, and, moreover that $(\lam,G^\lam)\cong(\rho_0,G^{\rho_0})$ is the
one-point compactification.

{\bf (i)}  If $G$ is abelian, then
neither of the converse quasi-containments to (\ref{eq:rachmancont}) hold;
see \cite[7.4.1]{grahamm}, for example.
This shows that the converse of (1) fails.  

{\bf (ii)} If $G=\mathrm{SL}_2(\Ree)$, then for any non-trivial complementary series
representation $\kappa_s$ ($0<s<1$), we have $\kappa_s\qc\rho_0$ by \cite[V.2.0.3]{howet},
for example, so
$\kappa\ec\rho_0\cong_e\lam$.  However it follows from Harish-Chandra's
trace formula \cite{harishchandra} that $\kappa_s\not\leq_w\lam$ for any $s$.
Hence the converse to (2) fails; in particular, Eberlein containment can hold
in a situation where weak containment fails.

Now, consider the representation $\kappa=\bigoplus_{0<s<1}\kappa_s$.
As indicated above, $\kappa\ec\lam$, and, since $\lam\cong\bar{\lam}$
it follows that $\kappa\we\lam$.  However it follows from \cite{milicic}
(see \cite[\S 7.6]{folland}) that $1\leq_w\kappa$, so $1\we\kappa$.  
However, $G$ is non-amenable so $1\not\leq_w \lam$
(we may also use Harish-Chandra's trace formula \cite{harishchandra} to see this).
Again by \cite[V.2.0.3]{howet}, this implies that $\wcoef_\lam\subset\fC_0(G)$,
and it follows that $1\not\leq_{we}\lam$.
\end{example}

We remark that by \cite[Theo. 6.1]{howem}, for any algebraic group $G$ over
a local field and any non-trivial irreducible representation $\pi$ we have
that $\pi\ec\rho_0$.

Recall that the definition of an $(\eta,H)$-regular compactification was given
before Proposition \ref{prop:ghahlaulos}.

\begin{theorem}\label{theo:regular}
A right topological compactification $(\del,S)$ of $G$ is such that $S$ has
an open group of units, if and only if $(\del,S)$ is a regular $(\eta,H)$-compactification
for some locally compact right topological group $H$ with proper right translations.
In particular,  if $\fC(S)\comp\del\supset\fC_0(G)$, then $\del(G)\cong G$ is the open
group of units in $S$.
\end{theorem}

\proof If the group $H$ of units of $S$ is open, then $H$ itself
is a locally compact right toplogical group.
Hence $(\del,S)$ is $(\del,H)$-regular.
If, on the other hand, $(\del,S)$ is $(\eta,H)$-regular, then
by (\ref{eq:gxdecomp}) in Proposition \ref{prop:ghahlaulos},
the group of units is isomorphic to $H$, and hence open.

By Proposition \ref{prop:ghahlaulos}, $(\del,S)$ is
$(\id, G)$-regular exactly when $\fC(S)\comp\del\supset\fC_0(G)$. \endpf

We note that not every Eberlein compactification has an open group
of units.  See Example \ref{ex:unitsnotopen} (iii).  Fortunately
$G^\fE$ has an open group of units.

\begin{corollary}\label{cor:regular1}
Let $G$ and $H$ be locally compact groups, $\pi\in\Sig_G$
satisfy $\lam_G\ec\pi$, and $\sig\in\Sig_H$ satisfy
$\lam_H\ec\sig$.  Then if $G^\pi\cong H^\sig$ as semitopological
semigroups, we must have $G\cong H$.  In particular,
$G^\fE\cong H^\fE$ if and only if $G\cong H$.
\end{corollary}

\proof  In light of Proposition \ref{prop:eberleincont}, our assumptions imply that
$\fC_0(G)=\fE(\lam_G)\subset\fE(\pi)+\Cee 1$.
Thus $\fC_0(G)\subset\fE(\pi)$, and the theorem above tells
us that the group of units of
$G^\pi$ is isomorphic to $G$.  The result follows.
\endpf

We observe that we can repeat the arguments above to see that
if $q:G\to G/\ker\pi$ is the quotient map, and
$\lam_{G/\ker\pi}\comp q\ec\pi$, then the group of
units of $G^\pi$ is isomorphic to $G/\ker\pi$.


\subsection{Subalgebras generated by matrix coefficients}\label{ssec:subalggenmatcoef}
We pass to the case of the algebras $\fA(\pi)$ shortly, with the aid
of a straightforward generalisation of \cite[Theo.\ 1]{walter1}.  Below we let
$\vn(\pi)=\vn_{\tau_\pi}$.  Also, if $\sig\qc\pi$,
we let $\sig^\pi:\vn_\pi\to\vn_\sig$
be the normal $*$-representation given by taking the adjoint of the inclusion map
$\ccoef_\sig\hookrightarrow\ccoef_\pi$; see, for example,
Proposition \ref{prop:xtovn}, above.  Note that this notation satisfies
$\sig^\pi(\pi(s))=\sig(s)$ for $s\iin G$.

\begin{theorem}\label{theo:walter}
Let $\pi\in\Sig_G$ and $x\in\vn(\pi)\cong\ccoef(\pi)^*$.
Then the following are equivalent

{\bf (i)}  $x\in\Phi_{\ccoef(\pi)}\cup\{0\}$;

{\bf (ii)} $(\sig\otimes\rho)^{\tau_\pi}(x)=\sig^{\tau_\pi}(x)\otimes\rho^{\tau_\pi}(x)$
for any $\sig,\rho\mq\pi$; and

{\bf (iii)} $(\pi\otimes\pi)^{\tau_\pi}(x)=\pi^{\tau_\pi}(x)\otimes\pi^{\tau_\pi}(x)$.

\noindent
Thus $\Phi_{\ccoef(\pi)}\cup\{0\}$ is a $*$-semigroup in $\vn(\pi)$, which
may be isomorphically identified with the $*$-semigroup
$\pi^{\tau_\pi}(\Phi_{\ccoef(\pi)})\cup\{0\}$ in $\vn_\pi$.
Also, if $x\in\pi^{\tau_\pi}(\Phi_{\ccoef(\pi)})$ has polar decomposition
$x=v|x|$, then $v,|x|\in\pi^{\tau_\pi}(\Phi_{\ccoef(\pi)})$ too.

If $\pi\otimes\pi\qc\pi$, i.e.\ $\tau_\pi\cong_q\pi$, then for $x\in\vn_\pi$
we have
\begin{equation}\label{eq:walter}
x\in\Phi_{\ccoef_\pi}\cup\{0\}\;\text{ if and only if }\;
(\pi\otimes\pi)^\pi(x)=x\otimes x.
\end{equation}
\end{theorem}

We shall indicate in Example \ref{ex:sxnotsemigroup}, below, that
$\Phi_{\ccoef(\pi)}$ is not itself a semigroup, in general.

\medskip
\proof  We first note that for typical elements of $\ccoef(\pi)$, say
$\inprod{\sig(\cdot)\xi}{\eta}$ and $\inprod{\rho(\cdot)\vartheta}{\zeta}$
where $\sig,\rho\mq\pi$, we have
$\inprod{\sig\otimes\rho(\cdot)\xi\otimes\vartheta}{\eta\otimes\zeta}
=\inprod{\sig(\cdot)\xi}{\eta}\inprod{\rho(\cdot)\vartheta}{\zeta}$.
Thus we have that $x\in\Phi_{\ccoef(\pi)}\cup\{0\}$ if and only if
\[
\inprod{(\sig\otimes\rho)^{\tau_\pi}(x)\xi\otimes\vartheta}{\eta\otimes\zeta}
=\inprod{\sig^{\tau_\pi}(x)\xi}{\eta}\inprod{\rho^{\tau_\pi}(x)\vartheta}{\zeta}
\]
where the latter expression is simply
$\inprod{\sig^{\tau_\pi}(x)\otimes\rho^{\tau_\pi}(x)\xi\otimes\vartheta}
{\eta\otimes\zeta}$.
Hence we have the equivalence of (i) and (ii).
That (ii) implies (iii) is clear.  That (iii) implies (i) follows from the computation above,
since $F_\pi$ generates $\ccoef(\pi)$. By the same fact we see
that $\pi^{\tau_\pi}$ must be injective on $\Phi_{\ccoef(\pi)}\cup\{0\}$.
The closure of $\Phi_{\ccoef(\pi)}$ under polar decomposition follows exactly
as in \cite[Theo.\ 1]{walter}, hence this happens in $\pi^{\tau_\pi}(\Phi_{\ccoef(\pi)})$ too.

We note that $\pi\otimes\pi\qc\pi$ implies that $\ccoef_{\pi}$ is itself an algebra,
in which case $\tau_\pi\qc\pi$, and hence $\tau_\pi\cong_q\pi$.
In this case $\pi^{\tau_\pi}$ is an isomorphism
with inverse $(\tau_\pi)^\pi$, and (\ref{eq:walter}) is immediate
from the equivalence of (i) and (ii) applied to elements $(\tau_\pi)^\pi(x)$.
\endpf

If $\sig\se\pi$,
Proposition \ref{prop:xtovn} provides a $*$-homomorphism
$\sig_e^\pi:\fE_\pi^*\to\vn_\sig$.
We note, moreover, that if $\sig\ec\pi$ then
$\sig_e^{\rho_\pi}:\fE(\pi)^*\cong\meas(G^\pi\setdif\{0\})\to\vn_\sig$
may be realised by the integral formula
\[
\inprod{\sig_e^\pi(\mu)\xi}{\eta}=\int_{G^\pi\setdif\{0\}}\inprod{\theta(x)\xi}{\eta}d\mu(x)
\]
where $\theta:G^\pi\to G^\sig$ is the factor map.
Also, if $\sig,\rho\me\pi$, then $\sig\otimes\rho\me\pi$, and hence
we can define $(\sig\otimes\rho)_e^{\tau_\pi}:\fA(\pi)^*\to
\vn_{\sig\otimes\rho}$.  Finally, since $\ccoef_\pi$ is dense in $\fE_\pi$,
$\pi_e^\pi:\fE_\pi^*\to\vn_\pi$ is injective, and hence we can identify
$\fE_\pi^*$ as a linear subspace of $\vn_\pi$

\begin{corollary}\label{cor:walter1}
Let $\pi\in\Sig_G$ and $\mu\in\fA(\pi)^*$.  Then the
following are equivalent

{\bf (i)} $\mu\in\Phi_{\fA(\pi)}\cup\{0\}$

{\bf (ii)} $(\sig\otimes\rho)_e^{\tau_\pi}(\mu)
=\sig_e^{\tau_\pi}(\mu)\otimes\rho_e^{\tau_\pi}(\mu)$
for any $\sig,\rho\me\pi$; and

{\bf (iii)} $(\pi\otimes\pi)_e^{\tau_\pi}(\mu)
=\pi_e^{\tau_\pi}(\mu)\otimes\pi_e^{\tau_\pi}(\mu)$.

\noindent Moreover,
$\pi_e^{\tau_\pi}:\Phi_{\fA(\pi)}\cup\{0\}\to\pi_e^{\tau_\pi}(\Phi_{\fA(\pi)}\cup\{0\})$
is a $*$-isomorphism and homeomorphism, in particular
latter set is a $*$-semigroup in $\vn_\pi$ which contains $G^\pi$.

If $\pi\cong_{se}\bar{\pi}$ we can replace (i) by

{\bf (i')} $\mu\in G^{\fA(\pi)}$ (in which case $\pi_e^{\tau_\pi}(\mu)=x$
for some $x\iin G^\pi\setdif\{0\}$, or $\mu=0$)

\noindent and $\ec$ can replace $\me$ in (ii).
\end{corollary}

\proof Since $\ccoef(\pi)$ is dense in $\fA(\pi)$, the computations
of Theorem \ref{theo:walter} can be repeated.
That $\Phi_{\fA(\pi)}\cup\{0\}$ is a semigroup on which $\pi_e^{\tau_\pi}$
is an isomorphism follows from the facts
above, or alternatively, Theorem \ref{theo:algspec} (i).
If $\pi\cong_{se}\bar{\pi}$,
then (i') follows from the standard fact that the spectrum of
$\fA(\pi)=\fE(\pi)\cong\fC_0(G^\pi\setdif\{0\})\cong G^\pi\setdif\{0\}$,
and, of course, $0$ is allowed in the case that $1\not\leq_e\pi$.
\endpf

We gain an augmentation of Corollary \ref{cor:regular1}.

\begin{proposition}\label{prop:walter0}
{\bf (i)} If $\pi\in\Sig_G$ satisfies $\lam\mq\pi$, then the semigroup
$\Phi_{\ccoef(\pi)}\cup\{0\}$ has an open group of units, isomorphic
to $G$.

{\bf (ii)} Suppose $G$ and $H$ are locally compact groups and
$\pi\iin\Sig_G$ and $\sig\in\Sig_H$ are such that $\lam_G\mq\pi$
and $\lam_H\mq\sig$.  If the semigroups $\Phi_{\ccoef(\pi)}\cup\{0\}$
and $\Phi_{\ccoef(\sig)}\cup\{0\}$ are continuously isomorphic, then $G\cong H$.
\end{proposition}

\proof We first prove (i).  We first observe that $\eps_{\ccoef(\pi)}(G)$ is open
in $\Phi_{\ccoef(\pi)}$.  Indeed, $\eps_{\ccoef(\pi)}:G\to\Phi_{\ccoef(\pi)}$
is continuous, and injective, since $\ccoef(\pi)\supset\falg$.  We note that
if $x\in\Phi_{\ccoef(\pi)}$ and $\langle u,x\rangle\not=0$ for some $u\iin\falg$,
then $x\in\eps_{\ccoef(\pi)}(G)$; indeed, there is $s\iin G$ for which
$\langle w,x\rangle=w(s)$ for $w\iin\falg$ by \cite[(3.34)]{eymard}, and since $\falg$
is an ideal in $\ccoef(\pi)$ --- i.e.\ an ideal in $\fsalg$ ---
it follows that $x=\eps_{\ccoef(\pi)}(s)$.  Now, if $s\in G$,
by \cite[(3.2)]{eymard}, find compactly supported $v\iin\falg$ so that $v(s)=1$.  If
$(x_i)\subset\Phi_{\ccoef(\pi)}$ is a net converging to $\eps_{\ccoef(\pi)}(s)$ then
$\lim_i\langle v,x_i\rangle=v(s)=1$, so eventually $\langle v,x_i\rangle
\not=0$.  Then, for such $i$, $x_i=\eps_{\ccoef(\pi)}(s_i)$ for some
$s_i\iin U=\{s\in G:v(s)\not=0\}$. Thus $\eps_{\ccoef(\pi)}(U)$ is
a neighbourhood of $\eps_{\ccoef(\pi)}(s)$, from which it follows that
$\eps_{\ccoef(\pi)}(G)$ is open in $\Phi_{\ccoef(\pi)}$.

Next we note that a simple modification of the proof of (\ref{eq:gxdecomp}) in
Proposition \ref{prop:ghahlaulos}, along with the characterisation of
$\eps_{\ccoef(\pi)}(G)$ in $\Phi_{\ccoef(\pi)}$, above gives that
$\Phi_{\ccoef(\pi)}\setdif\eps_{\ccoef(\pi)}(G)$ is an
ideal in $\Phi_{\ccoef(\pi)}\cup\{0\}$.
Thus we see that $\eps_{\ccoef(\pi)}(G)\cong G$ is the subgroup
of the identity $\tau_\pi(e)$, in $\Phi_{\ccoef(\pi)}\cup\{0\}$, and is open.

The result (ii) follows immediately from (i). \endpf

It is not generally true that $\Phi_{\fA(\pi)}=\Phi_{\ccoef(\pi)}$.
We call $\Phi_{\ccoef(\pi)}\setdif\Phi_{\fA(\pi)}$ the {\it Wiener-Pitt
part} of $\Phi_{\ccoef(\pi)}$.  If $\Phi_{\ccoef(\pi)}=\Phi_{\fA(\pi)}$
we will say that $\pi$ is {\em spectrally natural}; if not we will say it is a
{\it Wiener-Pitt} representation.  This nomenclature is motivated by
the example of Wiener and Pitt, exhibiting a non-neagitive element of $\fsal{\Zee}$
whose spectrum caontains $i$; see \cite[VII 9.1]{katznelson}.
Whenever $G$ admits a non-compact abelian quotient $H\cong G/N$,
then $G$ admits Wiener-Pitt representations.  For example consider
$\ome_H\comp q$ or even $\rho_0\comp q$ ($\rho_0$ is the Rajchman
representation of $H$ defined in Example \ref{ex:rachman}) where
$q:G\to H$ is the quotient map; see \cite[6.4.1]{rudinB} or \cite[8.2.3]{grahamm}
for example.  On the other hand, if $G$ is a connected semi-simple Lie group with finite centre,
then it follows from \cite{cowling} that $\ome_G$ is spectally natural.

Part (ii) of the following, if true, would be useful to us in Section \ref{sec:examples}.

\begin{conjecture}\label{conj:nonwp}
{\bf (i)} If $\sig\qc\pi$ and $\pi$ is spectrally natural, then $\sig$ is spectrally natural.

{\bf (ii)} If $\sig\qc\lam$, where $\lam$ is the left regular representation,
then $\sig$ is spectrally natural.
\end{conjecture}

\noindent We observe that if $\sig\qc\lam$ then $\sig\cmq\lam$ as well.
It follows \cite[Theo.\ 2.1]{bekkals} that $\ecoef(\sig)\cong\ccoef(G/\ker \sig)$
and $\ker\sig$ is a compact subgroup of $G$.  Hence
$\ccoef(\sig)$ is a closed translation-invariant subalgebra of
$\ecoef(\sig)\cong\ccoef(G/\ker \sig)$,
and thus $\fA(\sig)$ is a closed translation invariant
subalgebra of $\fC_0(G/\ker\sig)$.  This may help in resolving (ii).

We end this section by relating for $\pi\in\Sig_G$, the
algebra $\fA(\pi)$ to the $\pi$-Eberlein compactification $G^\pi$.

\begin{theorem}\label{theo:boundary}
If $\pi\in\Sig_G$, then
$\partial_{\fA(\pi)}=G^{\fA(\pi)}\setdif\{0\}$.
\end{theorem}

\proof
We shall find it convenient to make the identification
$G^{\fA(\pi)}\cong G^{\tau_\pi}$, which may be facilitated by
applying Corollary \ref{cor:walter1} to $\tau_\pi$ instead of to
$\pi$ itself.
By Theorem \ref{theo:algspec} we have that
$G^{\tau_\pi}\supset\partial_{\fA(\pi)}$, and it suffices to show that
$\partial_{\fA(\pi)}\cap\tau_\pi(G)\not=\varnothing$
to see that $\partial_{\fA(\pi)}=G^{\tau_\pi}\setdif\{0\}$.

For $\xi\in\fH_{\tau_\pi}$ with $\norm{\xi}=1$ we let
\[
u_\xi(s)=\frac{1}{2}\inprod{\tau_\pi(s)\xi}{\xi}+
\frac{1}{2}\inprod{\tau_\pi(s)\xi}{\xi}^2
\]
so $u_\xi$ is a positive definite
element of $\fA(\pi)$ with $\norm{u_\xi}_\infty=u_\xi(e)=1$.
We have for $x\in G^{\tau_\pi}$ that $|\inprod{x\xi}{\xi}|\leq 1$ and hence
\[
1=|\langle u_\xi,x\rangle|=\frac{1}{2}|\inprod{x\xi}{\xi}+\inprod{x\xi}{\xi}^2|
\quad \Leftrightarrow \quad \inprod{x\xi}{\xi}=1.
\]
Hence $K_\xi=\{x\in G^{\tau_\pi}:\inprod{x\xi}{\xi}=1\}$
is the maximal norming set for $u_\xi$, i.e.\ $K_\xi$ is the set of elements
$x\iin G^{\tau_\pi}$ for which $|\langle u_\xi,x\rangle |=1$.

Now suppose $B\subset G^{\tau_\pi}$ is a boundary for $\fA(\pi)$.
Note that $B\cap K_\xi$ is compact for every $\xi$, and non-empty
by considerations above.   For each finite
collection $\xi_1,\dots,\xi_n$ of norm one vectors we have that
$\bigcap_{j=1}^n K_{\xi_j}$ is maximally norming for
$u_{\xi_1}+\dots+u_{\xi_n}$
and hence $\bigcap_{j=1}^n(B\cap K_{\xi_j})$ is non-empty.
By the finite intersection property
$\bigcap\{B\cap K_\xi:\xi\in\fH_{\tau_\pi},\norm{\xi}=1\}\not=\varnothing$.
However, $\bigcap\{K_\xi:\xi\in\fH_{\tau_\pi},\norm{\xi}=1\}=\{\tau_\pi(e)\}$,
thus $\tau_\pi(e)\in B$.  Hence $\tau_\pi(e)\in \partial_{\fA(\pi)}$.  \endpf

We remark that if $\lam\me\pi$, i.e.\ $\fC_0(G)\subset\fA(\pi)$, then
each $\pi(s)$, $s\iin G$, is a strong boundary point, and the theorem
above is trivial.

\subsection{Eberlein weak containment and amenability}\label{ssec:ebweakcontamen}
Let us briefly review the position of almost periodic functions within
$\fE(G)$.  This approach is a variant of the one taken in \cite[\S 2]{rundes}.
Let $\fA\fP(G)$ denote the C*-algebra of almost periodic functions on $G$,
and $(\eps_{\fA\fP},G^{\fA\fP})$ denote the almost periodic compactification.
Let $\hat{G}_F$ denote the collection of finite dimensional irreducible
representations in $\Sig_G$, and $\pi_F=\bigoplus_{\sig\in\what{G}_F}\sig$.
Then $\ccoef_F(G)=\ccoef_{\pi_F}
=\ccoef(G^{\fA\fP})\comp\eps_{\fA\fP}=\fsalg\cap\fA\fP(G)$ (see
\cite[(2.20)]{eymard}), and hence $\fE(\pi_F)=\fA\fP(G)$.
A straightforward adaptation of \cite[Thm.\ 2.22]{burckel} shows that there 
is a minimal idempotent $e_{\fA\fP}$ in the semitopological
semigroup $G^\fE$ for which $\fA\fP(G)=\fC(G^{\fA\fP})\comp\eps_{\fA\fP}=e_{\fA\fP}\mult\fE(G)$.
Combining comments above we see that
$\ccoef_F(G)=\ccoef(G^{\fA\fP})\comp\eps_{\fA\fP}=e_{\fA\fP}\mult\fsalg$, qua subspaces of $\fE(G)$.
Moreover, $f\mapsto e_{\fA\fP}\mult f$
is a $*$-homomorphism since $e_{\fA\fP}$ is the minimal idempotent in $G^\fE=\Phi_{\fE(G)}$.
With the identifications in Corollary \ref{cor:walter1} applied to $\pi=\ome_G$,
we see that $e_{\fA\fP}\iin G^{\fA\fP}$ corresponds to the central projection
$p_F\iin\wstarg$ which covers $\pi_F$.  (The role of $p_F$
is discussed in \cite{walter}, where it is denoted $z_F$, but it is not discussed whether
this projection is in $G^\fE$.)  In particular we obtain closed ideals
\[
\ccoef_{PI}(G)=\{f\in\fsalg:p_F\mult f=0\}\aand
\fE_0(G)=\{f\in\fE(G):e_{\fA\fP}\mult f=0\}
\]
of $\fsalg$ and $\fE(G)$, respectively.  Here $\ccoef_{PI}(G)$ stands for
the ``purely infinite'' part of $\fsalg$; this is conjugation closed and translation invariant
since it is $\fsalg\cap\fE_0(G)$, and hence by \cite[(3.17)]{arsac} there
is a representation $\pi_{PI}$ for which $\ccoef_{PI}(G)=\ccoef_{\pi_{PI}}$.
By \cite[(3.12) \& (3.14)]{arsac},  $\pi\qc\pi_{PI}$ if and only if $\pi$ contains no
finite dimensional subrepresentations.
We note that $\fE_0(G)$ is exactly the space
of functions in $\fE(G)$ whose absolute values 
are in the kernel of the invariant mean on $\fW\fA\fP(G)$; this can be adapted from
\cite[Cor.\ 2.18]{burckel}.  These ideals give rise to
``semi-direct product'' decompositions
\[
\fsalg=\ccoef_{PI}(G)\oplus_{\ell^1}\ccoef_F(G)\aand
\fE(G)=\fE_0(G)\oplus\fA\fP(G).
\]

We may regard the following as
a refinement of Theorem \ref{theo:cstarspec} (iii).

\begin{theorem}\label{theo:algcontone}
If $\pi\in\Sig_G$, then the following are equivalent:

{\bf (i)} $\sig\qc\pi$ for some $\sig\iin\what{G}_F$;

{\bf (i')} $\sig\se\pi$ for some $\sig\iin\what{G}_F$;

{\bf (ii)} $1\cmq\pi$;

{\bf (ii')} $1\ec\pi$;

{\bf (iii)} $\fE_\pi\cap\fA\fP(G)\not=\{0\}$;

{\bf (iii')} $\fE(\pi)\cap\fA\fP(G)\not=\{0\}$; and

{\bf (iv)}  $0\not\in G^\pi$.
\end{theorem}

We remark that this holds for $\pi=\lam$ if and only if $G$ is compact.

\medskip
\proof That (i) implies (i'), and (ii) implies (ii') are obvious.  If (i) holds
then $1\qc\sig\otimes\bar{\sig}\cmq\pi$, so (ii) holds.  Similarly
(i') implies (ii').   
Condition (ii') implies that $1\in\fE(\pi)\cap\fA\fP(G)$, and hence we obtain (iii'). 
Theorem \ref{theo:cstarspec} (iii) gives the equivalence of (ii')
and (iv).  It remains to prove that (iii) implies (i), and (iii') implies (iii).


If it were the case that $p_F\mult\ccoef_\pi=\{0\}$, then for any
$u\iin\fE_\pi$ and any sequence $(u_n)\subset\ccoef_\pi$
for which $\lim_{n\to\infty}\unorm{u-u_n}=0$, we would have $e_{\fA\fP}\mult u
=\lim_{n\to\infty}e_{\fA\fP}\mult u_n=\lim_{n\to\infty}p_F\mult u_n=0$, thus
$e_{\fA\fP}\mult\fE_\pi=\{0\}$.
If (iii) holds then $e_{\fA\fP}\mult\fE_\pi=\fE_\pi\cap\fA\fP(G)\not=\{0\}$, and hence
$\ccoef_\pi\cap\ccoef_F(G)=p_F\mult\ccoef_\pi\not=\{0\}$, whence there is
$\rho\iin\Sig_G$ for which $\ccoef_\pi\cap\ccoef_F(G)=\ccoef_\rho$.
Since every $\rho$ for which $\ccoef_\rho\subset\ccoef_F(G)$ is completely
reducible, we see that statement (i) follows from (iii).

If it were the case that (iii) did not hold, i.e.\
$\fE_\pi\cap\fA\fP(G)=\{0\}$, then we would have that
$\ccoef_\pi\cap\ccoef_F(G)=\{0\}$.  Hence by \cite[(3.12)]{arsac}, 
$\pi$ would contain no subrepresentation of $\pi_F$, in which case $\pi\qc\pi_{PI}$, and hence
$\ccoef_\pi\subset\ccoef_{PI}(G)$.
Since $\ccoef_{PI}(G)$
is a closed and conjugation-closed subalgebra, in fact ideal, of $\fsalg$, it would follow
that $\ccoef(\pi)\cap\ccoef_F(G)=\{0\}$.  In this case, we would have
by the same argument of the the above paragraph
that $e_{\fA\fP}\mult\fE(\pi)=\{0\}$, hence we would obtain $\fE(\pi)\cap\fA\fP(G)=\{0\}$, which 
would violate (iii').  Thus (iii') implies (iii).
\endpf

As a consequence of the above theorem,
we observe that for a finite dimensional, hence totally reducible, representation $\pi$, 
and any $\sig$ in $\Sig_G$, that $\pi\qc\sig$ if and only if $\pi\se\sig$.  Hence the diagram
(\ref{eq:diagram}) may be simplified by noting that each implication between
the first and second rows is an equivalence.  

Following \cite[Theo.\ 5.1]{bekka}, we say that
$\pi\iin\Sig_G$ is {\it amenable} if and only if $1\wc\pi\otimes\bar{\pi}$.

\begin{theorem}\label{theo:amenrep}
Let $\pi\in\Sig_G$.  Then the following are equivalent:

{\bf (i)} $1\we\pi$;

{\bf (ii)} $\sig\we\pi$ for some $\sig\iin\what{G}_F$; and

{\bf (iii)} $\fE\fB(\pi)\cap\fA\fP(G)\not=\{0\}$.

\noindent Moreover, each of the above conditions implies each of the following
conditions

{\bf (iv)} $\pi$ is amenable; and

{\bf (v)} $\fE\fB(\pi\otimes\bar{\pi})\cap\fA\fP(G)\not=\{0\}$;


\noindent which are equivalent to one another;
and each of conditions {\rm (i)-(iii)} is implied by {\rm (iv)} or {\rm (v)},
provided that $\pi\otimes\bar{\pi}\leq_w \pi$
\end{theorem}

\proof   
First, that (i) implies (ii) is clear.
If (ii) holds, then $1\leq\sig\otimes\bar{\sig}$, so
$1\in F_{\sig\otimes\bar{\sig}}\subset\alg(F_\sig+F_{\bar{\sig}})$,
whence we get (iii).

Fell continuity of conjugation and (\ref{eq:inversioninvariance}) tell us that
$\wcoef_{\bar{\pi}}=\bar{\wcoef}_{\pi}$, i.e.\ $\ome_{\bar{\pi}}=\bar{\ome}_\pi$.
Hence condition (i) is the same as saying $1\ec\ome_\pi$, while (iii)
is the same as saying $\fE(\ome_\pi)\cap\fA\fP(G)\not=\{0\}$.  Hence the
equivalence of (i) and (iii) follows from the equivalence of (ii') and (iii') of 
Theorem \ref{theo:algcontone}, above.
Moreover, we see by the equivalence of (i) and (ii') of Theorem \ref{theo:algcontone},
that condition (i) of the present theorem 
implies that $\sig\qc\ome_\pi$ for some $\sig$ in $\what{G}_F$.
Hence $1\leq \sig\otimes\bar{\sig}\qc\ome_\pi\otimes\ome_{\bar{\pi}}
\qc\ome_{\pi\otimes\bar{\pi}}$.  The last containment holds since
$\wcoef_\pi \wcoef_{\bar{\pi}}\subset \wcoef_{\pi\otimes\bar{\pi}}$ by virtue of the facts
that $F_\pi F_{\bar{\pi}}= F_{\pi\otimes\bar{\pi}}$ and
multiplication is weak*-continuous on $\fsalg$.  Hence we see (i) implies (iv).



If (v) holds, then by the equivalence of (i) and (iii') of Theorem \ref{theo:algcontone},
$\sig\qc\ome_{\pi\otimes\bar{\pi}}$ for some $\sig$ in $\what{G}_F$.  Hence it follows that
$\pi\otimes\bar{\pi}$ is amenable, and therefore $\pi$ is amenable (see
\cite[Theo.\ 1.3]{bekka} and \cite[Prop.\ 2.7]{stokke0}).  Hence we get (iv).  Conversely, if
(iv) holds, then $1\in\wcoef_{\pi\otimes\bar{\pi}}\subset
\fE\fB(\pi\otimes\bar{\pi})\cap\fA\fP(G)$ so (v) holds.


If we assume that $\pi\otimes\bar{\pi}\leq_w \pi$, and (iv) holds, then
$1\leq_w\pi\otimes\bar{\pi}\leq_w \pi$.   Thus (i) holds. 
\endpf

We remark that there are no evident containment relations between
$\fE\fB(\pi)$ and $\fE\fB(\pi\otimes\bar{\pi})$, in general.  However, if  $1\leq_w\pi$ we have
$\fE\fB(\pi)\subset\fE\fB(\pi\otimes\bar{\pi})$.


We now define the {\it reduced Eberlein algebra} and {\it reduced
Eberlein compactification} by
\[
\fE_r(G)=\wbar{\wcoef_r(G)}^{\unorm{\cdot}}\aand
(\eps_{\fE_r},G^{\fE_r})=(\eps_{\fE_r(G)},G^{\fE_r(G)})
\]
where $\wcoef_r(G)=\wcoef_\lam$ and
$\lam:G\to\fU(\bl^2(G))$ is the left regular representation. 
The following is an analogue  of Hulanicki's theorem \cite{hulanicki} ---
{\em $G$ is amenable if and only if each continuous positive definite
function can be approximated, uniformly on compacta, by positive definite
functions associated with the left regular representation} ---
which is equivalent to saying that $\ome_G\wc\lam$, or that $1\wc\lam$. 
However, parts (viii), (ix) and (x) have no obvious analogue in that context.

\begin{theorem}\label{theo:amenability}
The following are equivalent:

{\bf (i)}  $G$ is amenable;

{\bf (ii)} $\fE_r(G)=\fE(G)$;

{\bf (iii)} $1\in\fE_r(G)$;

{\bf (iv)} $\meas(G^{\fE_r}\setdif\{0\})\cong\fE_r(G)^*$ and
$\meas(G^\fE)\cong\fE(G)^*$ are isomorphic dual Banach algebras;

{\bf (v)} $\meas(G^{\fE_r}\setdif\{0\})\cong\fE_r(G)^*$ admits a weak* continuous character;

{\bf (vi)} $\sig\we\lam$ for some $\sig\in\hat{G}_F$; and

{\bf (vii)} $\pi\we\lam$ for every $\pi\iin\Sig_G$.

\noindent
If we further assume that $\hat{G}_F\setdif\{1\}\not=\varnothing$, then (i)-(vii), above,
are equivalent to

{\bf (viii)} $(\sig,G^\sig)\leq(\eps_{\fE_r},G^{\fE_r})$ for some $\sig\iin\hat{G}_F\setdif\{1\}$;

{\bf (ix)} $(\pi,G^\pi)\leq(\eps_{\fE_r},G^{\fE_r})$ for every $\pi\iin\Sig_G$; and

{\bf (x)} $(\eps_{\fE_r},G^{\fE_r})\cong(\eps_\fE,G^\fE)$.
\end{theorem}

\proof If (i) holds, then $\wcoef_r(G)=\fsalg$ and hence (ii) holds.  That (ii) implies
both (iii) and (iv) is obvious.
Condition (iii) is that $1\we\lam$, which by Theorem \ref{theo:amenrep}
implies that $\lam$ is amenable.  Hence by \cite[Theo.\ 2.2]{bekka}
(see also \cite[p.\ 260]{fell}, since we have used an equivalent definition
of ``amenable representation'') we obtain (i).
Since $\pi\we\ome_G$ holds generally, and (ii) says that $\ome_G\cong_{we}\lam$,
(ii) and (vii) are immediately equivalent.

If (iii) holds then, by Theorem \ref{theo:cstarspec}, $0\not\in G^{\fE_r}$.  Moreover
$\mu\con\nu(1)=\mu(\nu\mult 1)=\mu(\nu(1)1)
=\mu(1)\nu(1)$ for $\mu,\nu\iin\meas(G^{\fE_r})\cong\fE_r(G)^*$,
so (v) holds.  Similarly, $\meas(G^\fE)\cong\fE(G)^*$ always admits
a weak* continuous character, so if (iv) holds, then so too must (v).
If (v) holds, then there is $h\iin\fE_r(G)$
such that $\mu\con\nu(h)=\mu(h)\nu(h)$ for each $\mu,\nu\iin
\meas(G^{\fE_r}\setdif\{0\})\cong\fE_r(G)^*$.  Thus we have
$h(st)=\del_{st}(h)=\del_s\con\del_t(h)=\del_s(h)\del_t(h)=h(s)h(t)$ for $s,t\iin G$.
Hence $h$ is a norm one character, so $h\in\what{G}_F$
with $h\we\lam$, and we have (vi).  Converesly, if (vi) holds then
by Theorem \ref{theo:amenrep} we obtain $1\we\lam$ and hence
we obtain (iii).

Condition (vii) implies condition (ix), which in turn implies condition (viii).
If there is $\sig$ in $\hat{G}_F\setdif\{1\}$ satisfying condition (viii),
then by virtue of Proposition \ref{prop:eberleincont} we have
that $\Cee 1\not= \fE_\sig\subset\fE_r(G)\oplus \Cee 1$, from which it follows that
$\fA\fP(G)\cap\fE_r(G)\supset \fE_\sig\cap\fE_r(G)\not=\{0\}$.
Hence from the equivalence of (iii) and (ii') in Theorem \ref{theo:algcontone},
$1\ec \ome_\lam$ which in turn gives (iii) in the present theorem.
The choice $\ome=\pi$ in (ix) gives that $(\eps_\fE,G^\fE)\leq
(\eps_{\fE_r},G^{\fE_r})$; whereas the converse comparison holds by
Proposition \ref{prop:eberleincont}.  Thus (ix) implies (x).
That (x) implies (ix) follows from Proposition \ref{prop:eberleincont}.
\endpf

If $G=\mathrm{SL}_2(\Ree)$, then $\hat{G}_F=\{1\}$ and
$(\eps_\fE,G^\fE)=(\eps_{\fW\fA\fP},G^{\fW\fA\fP})$ is the one-point compactification;
see Example \ref{ex:rachman} (ii) and \cite{chou1}.
Hence $(\eps_{\fE_r},G^{\fE_r})=(\eps_{\fW\fA\fP},G^{\fW\fA\fP})$, in this case, and
the condition $(\pi,G^\pi)\leq(\eps_{\fE_r},G^{\fE_r})$ of (ix) is satisfied,
despite that $G$ is not amenable.

\subsection{Universal property of the Eberlein compactification}\label{ssec:universal}

We recall the definition of the (CH)-compactification from Section \ref{ssec:CH}.
The following theorem appears not to be new, but appears in a dual
form, {\it \`{a} la} Theorem \ref{theo:comparison}, in \cite[Thm.\ 3.12]{megrelishvili}.
Moreover no assumption of local compactness of the the topological group
$G$ is made in \cite{megrelishvili}.  Our proof appears to hold in that setting as well,
and is sufficiently different to merit inclusion.

\begin{theorem}\label{theo:universal}
The Eberlein compactification $(\eps_\fE,G^\fE)$ is the (CH) compactification
$(\eps_{\fC\fH},G^{\fC\fH})$, and hence it is universal amongst (CH) compactifications.
\end{theorem}

\proof We first note that (CH) compactifications are exactly those of the
form $(\pi,G^\pi)$ where $\pi\in\Sig_G$.  Indeed each such
$(\pi,G^\pi)$ is clearly a (CH) compactification.  Conversely, if
$\del:G\to\fB(\fH)_{\norm{\cdot}\leq 1}$ is a homomorphism, then
$p=\del(e)$ is a contractive idempotent, hence a projection.
For any $s\iin G$ let $u=\del(s)$ and $v=\del(s^{-1})$.  Then
$uv=vu=p$.  Moreover if $\xi\in p\fH$ then $\norm{\xi}=\norm{vu\xi}
\leq \norm{u\xi}\leq\norm{\xi}$, and if $\xi\in(1-p)\fH$ then
$u\xi=u(1-p)\xi=0$, which shows that $u$ is a partial isometry with
support and range projection $p$, and hence $v=u^*$.
We note that $p$ commutes with $S=\wbar{\del(G)}^{w^*}$.  We let
$\fH_\pi=p\fH$ and $\pi=p\del(\cdot)|_{p\fH}$,
so $\pi$ is a unitary representation.  Since $S=pSp$ we see that
$x\mapsto px|_{\fH_\pi}$ is a continuous bijection from $S$ onto $G^\pi$,
intertwining $\del$ and $\pi$.  Hence $(\del,S)\cong(\pi,G^\pi)$.

It is immediate from Proposition \ref{prop:eberleincont}, and the fact that
each $\pi\in\Sig_G$ satisfies $\pi\qc\ome_G$, hence $\pi\ec\ome_G$, that
each (CH) compactification is a factor of $(\eps_\fE,G^\fE)\cong(\ome_G,G^{\ome_G})$.
\endpf

Given $\pi\in\Sig_G$, we call a compactification $(\del,S)$
an (F$\pi$)-compactification (``factor of $\pi$'') if $(\del,S)\leq(\pi,G^\pi)$.
In keeping with Theorem \ref{theo:universal}, we have that an
(FCH)-compactification is an $(F\ome_G$)-compactification.
It would be interesting to know if every (FCH)-compactification
of $G$ is itself a (CH)-compactification.  It may be necessary
to condsider only involutive (FCH)-compactifications.



We can use Theorem \ref{theo:universal} to give examples of involutive
compact semitopological semigroups which cannot be faithfully represented
on Hilbert spaces.  We gain an extension of \cite[Theo.\ 4.7]{megrelishvili},
where it is observed via \cite{rudin} that the compact involutive semigroup $\Zee^w$
cannot be represented faithfully as contractions on a Hilbert space.
We say that $G$ is {\it totally minimal} if every continuous
homomorphism into another topological group has closed range.
It is shown in \cite[Thm.\ 2.5]{mayer} that connected totally minimal groups are precisely
those which are inductive limits of groups of the form $R\rtimes N$, where
$R$ is a reductive Lie group acting on a nilpotent Lie group $N$ with no
non-trivial fixed points.

\begin{corollary}\label{cor:universal2}
Suppose that $G$ is either (a) nilpotent, (b) has an inner automorphism invariant
compact neighbourhood of the identity, i.e.\ is an [IN]-group, or (c) is connected but not
totally minimal.  Then $(\eps_{\fW\fA\fP},G^{\fW\fA\fP})\not\leq (\eps_\fE,G^\fE)$.  In particular
$(G^{\fW\fA\fP})^{\fC\fH}$ is a proper quotient of $G^{\fW\fA\fP}$, in this case.
\end{corollary}

\proof For cases (a) and (b)
it is shown in \cite{chou}, and for case (c) it is shown in \cite[Thm.\ 4.5]{mayer},
that $\fW\fA\fP(G)\not\subset\fE(G)$.
Hence that $(\eps_{\fW\fA\fP},G^{\fW\fA\fP})\not\leq (\eps_\fE,G^\fE)$ is immediate
from Theorem \ref{theo:comparison} (ii).  If it were the case that there
were an isomorphism $\theta:G^{\fW\fA\fP}\to (G^{\fW\fA\fP})^{\fC\fH}$, then
$(\theta\comp\eps_{\fW\fA\fP},(G^{\fW\fA\fP})^{\fC\fH})$
would be a (CH)-compactification of $G$, which would imply that
$(\eps_{\fW\fA\fP},G^{\fW\fA\fP})\leq (\eps_\fE,G^\fE)$, which contradicts results above.  \endpf

It is shown in \cite[Thm.\ 4.5]{mayer} that
all connected groups for which $(\eps_\fE,G^\fE)\cong(\eps_{\fW\fA\fP},G^{\fW\fA\fP})$
have Eberlein compactifications of the
form illustrated in Proposition \ref{prop:ilies}, below.

We also have a complementary
``minimality property'' of Eberlein type compactifications.

\begin{proposition}\label{prop:minimal}
Let $(\del,S)$ be a right topological compactification of $G$, and
$\pi\in\Sig_G$.  Then the following are equivalent:

{\bf (i)} there is a continuous homomorphism $\pi_S:S\to(\fB(\fH_\pi),w^*)$
such that $\pi_S\comp\del=\pi$;

{\bf (ii)} $(\pi,G^\pi)\leq(\del,S)$;

{\bf (iii)} $\fC(S)\comp\del\supset\ccoef_\pi$.

In particular, $(\eps_\fE,G^\fE)$ is the minimum compactification
for which (i), above, holds for all $\pi\iin\Sig_G$.
\end{proposition}

\proof  That (i) implies (ii) is by definition of the ordering of compactifications.
Note that since $\fB(\fH_\pi)_{\norm{\cdot}\leq 1}$ is weak* compact,
$\pi_S(S)\subset \fB(\fH_\pi)_{\norm{\cdot}\leq 1}$.  Since
$\fE(\ccoef_\pi)=\fE(\pi)$,  it follows
from Theorem \ref{theo:comparison} that (ii) implies (iii).
If (iii) holds, a straighforward modification of Proposition
\ref{prop:xtovn} provides a weak*-weak* continuous map
$\del_S:\fC(S)^*\cong(\fC(S)\comp\del)^*\to\vn_\pi$ for which
$\del_S\comp\eps_{\fC(S)}\comp\del=\pi$.
Let $\pi_S=\del_S\comp\eps_S$ and we obtain (i).

Thus if $(\del,S)$ satisfies (i) for all $\pi\iin\Sig_G$,
then with the choice of $\pi=\ome$ we obtain that
$(\del,S)\geq(\eps_\fE,G^\fE)$.  Of course, if
$(\del,S)\geq(\eps_\fE,G^\fE)$, then (i) holds
for all $\pi\iin\Sig_G$.  \endpf


\section{Examples}\label{sec:examples}

As in the previous section, we will always let $G$ denote a locally compact group,
unless otherwise indicated.

\subsection{Spine type compactifications}\label{ssec:spine}
Suppose $H$ is a locally compact group and $\eta:G\to H$ is a continuous
homomorphism with dense range.  We will call $(\eta,H)$ a {\it locally compact
completion} of $G$.  Two locally compact completions
$(\eta_j,H_j)$, $j=1,2$, have a {\it mutual quotient}
if there are compact normal subgroups $K_j\subset H_j$ for which
$H_1/K_1\cong H_2/K_2$, via a bicontinuous isomorphism $\theta$ for
which $\theta\comp q_i\comp\eta_i=q_j\comp\eta_j$, where $q_j:H_j\to H_j/K_j$
is the quotient map for each $j$.  We shall define the {\it subdirect product}
of locally compact completions $(\eta_j,H_j)$, $j=1,2$, as the pair
$(\eta_1\cross\eta_2,\wbar{\{(\eta_1(s),\eta_2(s)):s\in G\}})$ and denote it by
$(\eta_1,H_1)\vee(\eta_2,H_2)$.

The following is adapted from \cite{ilies}.  We let $\lam_H:H\to\bl^2(H)$ denote
the left regular representation of $H$.

\begin{proposition}\label{prop:ilies}
Let $(\eta_j,H_j)_{j\in J}$ be a family of locally compact completions of $G$, and
\[
\lam_J=\bigoplus_{j\in J}\lam_j\text{ where }\lam_j=\lam_{H_j}\comp\eta_j.
\]

{\bf (i)} No two distinct
$(\eta_i,H_i)$ and $(\eta_j,H_j)$ are mutually quotient if and only if
\[
\ccoef_{\lam_J}=\ell^1\text{-}\bigoplus_{j\in J}\ccoef(\lam_j).
\]
Moreover, each $\ccoef(\lam_j)=\ccoef_{\lam_j}=\ccoef(H_j)\comp\eta_j\cong\ccoef(H_j)$.

{\bf (ii)} If for each $i,j\iin J$, $(\eta_i,H_i)\vee(\eta_j,H_j)\cong
(\eta_k,H_k)$ for some $k\iin J$, we write $k=i\vee j$ and
$(J,\vee)$ is a semilattice.  Further
$\ccoef_{\lam_J}=\ccoef(\lam_J)$ is an algebra of functions, graded over $(J,\vee)$ in the sense
that $\ccoef(\lam_i)\ccoef(\lam_j)\subset\ccoef(\lam_{i\vee j})$.
Moreover, $i\vee j=j$ --- we write $i\leq j$ ---
if and only if  there is a continuous homomorphism
$\eta_i^j:H_j\to H_i$ such that $\eta_i^j\comp\eta_j=\eta_i$.

{\bf (iii)} Suppose both (i) and (ii) hold and also that the semilattice
$(J,\vee)$ is complete. Then for any $i,j\iin J$ either

\parbox[t]{4.5in}{{\bf (a)} there
is an element $i\wedge j\iin J$ which satisfies $i\wedge j\leq i$,
$i\wedge j\leq j$, and for any $k\iin J$ for which $k\leq i$ and
$k\leq j$, then $k\leq i\wedge j$ too; or }

\parbox[t]{4.5in}{{\bf (b)} there is no $k\iin J$ for which $k\leq j$ and $k\leq i$.}

\noindent If (a) always holds for all $i,j$ then $(J,\wedge)$ is a semilattice; otherwise
we can adjoin an identity $o$, so $o\vee j=j$ for each $j\iin J$, and $(J\cup\{o\},\wedge)$
is a semilattice.  We obtain identifications
\[
\Phi_{\ccoef(\lam_J)}\cong G^{\lam_J}\setdif\{0\}\cong\bigsqcup_{j\in J} H_j.
\]
The semigroup structure on $G^{\lam_J}$ is given for $s_i\iin H_i$ and $s_j\iin H_j$ by
\begin{equation}\label{eq:semigroupstructure}
s_is_j=\begin{cases} \eta^i_{i\wedge j}(s_i)\eta^j_{i\wedge j}(s_j) &\text{ if (a) holds} \\
0&\text{ if (b) holds.}\end{cases}
\end{equation}
The topology is given by basic open neighbourhoods of a point $s\iin H_j$ that are of
the form
\[
V_j\sqcup\left\{t\in \bigsqcup_{i\in J,i>j}H_i:\eta^i_{i_k}(t)\in W_{i_k}\text{ if }
t\in H_i\text{ for }i\geq i_k\right\}
\]
where $V_j$ is an open neighbourhood of $s$ in $H_j$, each $i_k>j$,
and each $W_{i_k}$ is a cocompact set in $H_{i_k}$.
\end{proposition}

\proof  Part (i) can proved as in \cite[\S 3.3]{ilies}.  Part (ii)
can be proved as in \cite[\S 3.2]{ilies}.  Part (iii) is proved as in
\cite[\S 4.1-4.3]{ilies}.   We remark on some details given there.
We let
\[
\mathfrak{HD}(J)=\left\{S\subset J:
\begin{matrix}
S\text{ is hereditary: }j\in S,i\leq j &\Rightarrow  & i\in S;\\
\aand S\text{ is directed: }i,j\in S &\Rightarrow &i\vee j\in S
\end{matrix}\right\}.
\]
For each $s\in\mathfrak{HD}(J)$ we obtain an inverse mapping system
$\{H_j,\eta^j_i:j\in S, i\leq j\iin S\}$ which gives rise to a projective limit
\[
H_S=\underset{j\in S}{\underleftarrow{\lim}}H_j
=\left\{(s_j)_{j\in S}\in\prod_{j\in S}H_j:\eta^j_i(s_j)=s_i\text{ if } i\leq j\iin S\right\}.
\]
Then  the proof of \cite[Theo.\ 4.1]{ilies} and remarks on \cite[p.\ 285]{ilies} gives
the structure of $G^{\lam_J}$and semigroup product by
\[
G^{\lam_J}\cong\bigsqcup_{S\in \mathfrak{HD}(J)}G_S,\quad
(s_j)_{j\in S_1}(t_j)_{j\in S_2}=
\begin{cases}(s_jt_j)_{j\in S_1\cap S_2}&\text{if }S_1\cap S_2\not=\varnothing\\
0&\text{if }S_1\cap S_2=\varnothing.\end{cases}
\]
The assumption that $(J,\vee)$ is complete allows that each
element of $\mathfrak{HD}(J)$ is principal, i.e.\ of the form $S_j=\{i\in J:i\leq j\}$,
in which case $H_{S_j}\cong H_j$.  In the event that $S_i\cap S_j\not=\varnothing$
for each $i,j$, we obtain (a); otherwise (b) holds.
In the case that (b) holds, for some $i,j$, we adjoin $o$ to $J$ so
$o\leq j$ for each $j\iin J$.  We let $\eta^j_o(s)=0$, and we obtain the product
(\ref{eq:semigroupstructure}) as in \cite[(4.8)]{ilies}.

The description of
the topology is immediate from the corrigendum \cite{iliesC} to
\cite[Theo.\ 4.2]{ilies}.  Translated into the
present terminology a basic open neighbourhood of
$s=(s_j)_{j\in S_0}\in H_{S_0}\subset G^{\lam_J}$
is of the form
\begin{equation}\label{eq:neighbourhood}
\left\{t\in G^{\lam_J}:
\begin{matrix}t=(t_i)_{i\in S}\in H_S\text{ for some }S\supseteq S_j\iin\mathfrak{HD}(J) \\
\text{where } t_j\in V_j\aand  t_{i_k}\in W_{i_k}\text{ if }i_k\in S\end{matrix}\right\}
\end{equation}
where $j\in S_0$, $V_j$ is an open neighbourhood of $s_j$, $i_1,\dots,i_k>j$ and
$W_{i_k}$ is a cocompact set in $H_{i_k}$ for each $k$.
\endpf

We note that since $G$ is dense in $G^{\lam_J}$, $\lam_J$ is spectrally natural
in the sense that $\Phi_{\fA(\lam_J)}=\Phi_{\ccoef(\lam_J)}$.  In fact,
it can be verified by way of the regularity condition for Fourier algebras
\cite[(3.2)]{eymard} that this algebra is regular on $G^{\lam_J}$.  We say that the
representation $\lam_J$ {\it spectrally regular} in this case.

\begin{example}\label{ex:sxnotsemigroup}
Here we show examples satisfying all of the assumptions of Proposition
\ref{prop:ilies} for which $G^{\lam_J}\setdif\{0\}=\Phi_{\fA(\lam_J)}$ is not a semigroup.

{\bf (i)} Let $H$ be any non-compact locally compact group and $G=H\cross H$.
Let $(\eta_o,H_o)=(\id,G)$, $H_l=H=H_r$, and $\eta_l(s,t)=s$, $\eta_r(s,t)=t$
for $(s,t)\in G$.  Then $J=\{o,l,r\}$, in the notation above, is the flat semilattice
given by $o\vee j=o$ for any $j\iin J$ and $r\vee l=o$.  Clearly there is no
$j\iin J$ for which $j\leq l$ and $j\leq r$, so $G^{\lam_J}\setdif\{0\}$ is not itself
a semigroup.

{\bf (ii)} Let $G=\Ree^n$.  Fix an inner product on $G$ and
for any subspace $L\subset G$, let $\eta_L$ be the orthogonal
projection onto $L$.  If  $\fL$ denotes the set of subspaces, then
$L_1\vee L_2=L_1+L_2$ and $L_1\wedge L_2=L_1\cap L_2$.
If we consider, for example, $\fL_k=\{L\in\fL:\dim L\geq k\}$,
for $k=0,\dots,n$, then $G^{\lam_{\fL_k}}\setdif\{0\}$ is a semigroup only if $k=0,n$.
\end{example}

\begin{example}\label{ex:unitsnotopen}
Let $\Que$ be the discrete rationals and
$G$ be the direct sum group $\Que^{\oplus\infty}$.
Let for each $n$, $\eta_n:G\to\Ree^n$ be the projection onto
the first $n$ coordinates.
We note that $(\eta_n,\Ree^n)\vee(\eta_m,\Ree^m)
=(\eta_{n\vee m},\Ree^{n\vee m})$ where $n\vee m=\max\{n,m\}$, and
if $m\leq n$ then $\eta^n_m:\Ree^n\to\Ree^m$ is the projection onto the
first $m$ coordinates.  In the notation of the proof of
Proposition \ref{prop:ilies}, the projective limit $H_\En$ is
isomorphic to the direct sum group $\Ree^{\oplus\infty}$
with inductive limit topology, hence is not locally compact.
We note that every element of $\mathfrak{HD}(\En)$
is principal except for $\En$ itself.  Thus we obtain the structure and
semigroup product
\[
G^{\lam_{\En}}\setdif\{0\}\cong\bigsqcup_{n\in\En\cup\{\infty\}}\Ree^{\oplus n},
\quad(s_j)_{j=1}^n(t_j)_{j=1}^m=(s_j+t_j)_{j=1}^{n\wedge m},
\,n,m\in\En\cup\{\infty\}.
\]
The basic open neighbourhoods
of $s=(s_j)_{j=1}^\infty\in\Ree^{\oplus\infty}$, as described in (\ref{eq:neighbourhood}), are
\[
\left\{t=(t_j)_{j=1}^m\in\bigsqcup_{m=n,\dots,\infty}\Ree^{\oplus m}:
(t_j)_{j=1}^n\in V_n\right\}
\]
where $V_n$ is an open neighbourhood of $(t_j)_{j=1}^n$ in $\Ree^{\oplus n}$.
The group of units is $\Ree^\infty$ and is plainly not open in $G^{\lam_{\En}}$.
\end{example}

We wish to extend Proposition \ref{prop:ilies} to include subrepresentations
of those of the type $\lam_J$.  

%

\begin{proposition}\label{prop:morespine}
Let $\{\pi_j,\fH_j\}_{j\in J}$ be a family of representations for which
there is a corresponding family $(\eta_j,H_j)$ of locally compact completions
of $G$ for which $\pi_j\qc\lam_{H_j}\comp\eta_j$ for each $j$.
Suppose that the assumptions of (i), (ii) and (iii) in Proposition \ref{prop:ilies}
hold and further that
\[
\ccoef_{\pi_J}=\ell^1\text{-}\bigoplus_{j\in J}\ccoef(\pi_j)
\]
is a subalgebra of $\ccoef_{\lam_J}$, where $\pi_J=\bigoplus_{j\in J}\tau_{\pi_j}$.
Then for each $i\leq j$ in $J$
there is a semigroup homomorphism $\kappa^j_i:\Phi_{\ccoef(\pi_j)}\cup\{0\}
\to \Phi_{\ccoef(\pi_i)}\cup\{0\}$ which satisfies $\kappa^j_i\comp\tau_{\pi_j}(s)
=\tau_{\pi_i}(s)$ for each $s\iin G$ and $\kappa^j_i(0)=0$.
Moreover, $\ccoef_{\pi_J}=\ccoef(\pi_J)$ and
\[
\Phi_{\ccoef(\pi_J)}\cong\bigsqcup_{j\in J}\Phi_{\ccoef(\pi_j)}.
\]
The semigroup structure on $\Phi_{\ccoef(\pi_J)}\cup\{0\}$ is given
for $x_i\iin \Phi_{\ccoef(\pi_i)}\cup\{0\}$ and $x_j\iin \Phi_{\ccoef(\pi_j)}\cup\{0\}$ by
\[
x_ix_j=\begin{cases} \kappa^i_{i\wedge j}(x_i)\kappa^j_{i\wedge j}(x_j)
&\text{ if (a) of Proposition \ref{prop:ilies} holds} \\
0&\text{ if (b) of Proposition \ref{prop:ilies} holds.}\end{cases}
\]
The topology is given as follows:  basic open neighbourhoods of a point
$x_j\iin \Phi_{\ccoef(\pi_j)}$ are of
the form
\[
V_j\sqcup\left\{x\in \bigsqcup_{i\in J,i>j}\Phi_{A(\pi_j)}:\kappa^i_{i_k}(x)\in W_{i_k}\text{ if }
x\in \Phi_{\ccoef(\pi_i)}\text{ for }i\geq i_k\right\}
\]
where $i_1,\dots,i_k> j$, each $W_{i_k}$ is a cocompact set in $\Phi_{\ccoef(\pi_{i_k})}$,
and $V_j$ is an open neighbourhood of $s_j$ in $\Phi_{\ccoef(\pi_j)}$.
\end{proposition}

\proof Whilst similar to Proposition \ref{prop:ilies}, the present result
cannot be deduced from the same proof from \cite{ilies}.  Thus we show some of the
details.

The assumption that $\ccoef_{\pi_J}$ is an algebra
and the assumptions of Proposition \ref{prop:ilies} (i) and (ii) imply that for
$i\leq j$ we have
\[
\ccoef(\pi_i)\ccoef(\pi_j)\subset \ccoef_{\pi_J}\cap \bigl(\ccoef(\lam_i)\ccoef(\lam_j)\bigr)
\subset\ccoef_{\pi_J}\cap \ccoef(\lam_j)=\ccoef(\pi_j).
\]
If $i\leq j$ and $x\in\Phi_{\ccoef(\pi_j)}$, define $\kappa^j_i(x)\iin\vn(\pi_i)$
by
\[
\langle u_i,\kappa^j_i(x)\rangle=\frac{\langle u_iu_j, x\rangle}{\langle u_j, x\rangle}
\]
for $u_i\iin\ccoef(\pi_i)$ and $u_j\iin\ccoef(\pi_j)$ for which $\langle u_j, x\rangle\not=0$.
This is independent of the choice of $u_j$, since if $u_j'\in\ccoef(\pi_j)$ then
\[
\langle u_iu_j, x\rangle\langle u_j', x\rangle
=\langle u_iu_ju_j', x\rangle
=\langle u_iu_j', x\rangle\langle u_j, x\rangle.
\]
We see that that $\kappa^j_i(x)\in\Phi_{\ccoef(\pi_i)}\cup\{0\}$:
if $u_i,u_i'\in\ccoef(\pi_i)$ and $u_j,u_j'$ are as above and with
$\langle u_j,x\rangle\langle u_j', x\rangle\not=0$, then
\[
\langle u_i,\kappa^j_i(x)\rangle\langle u_i',\kappa^j_i(x)\rangle
=\frac{\langle u_iu_j, x\rangle}{\langle u_j, x\rangle}
\frac{\langle u_i'u_j', x\rangle}{\langle u_j', x\rangle}
=\frac{\langle u_iu_i'u_ju_j', x\rangle}{\langle u_ju_j', x\rangle}
=\langle u_iu_i',\kappa^j_i(x)\rangle.
\]
It is obvious that $\kappa^j_i(0)=0$ and straightforward to
check that $\kappa^j_i\comp\eta_j=\eta_i$ on $G$.

We wish to show that $\kappa^j_i$ is multiplicative on $\Phi_{\ccoef(\pi_j)}\cup\{0\}$.
First,  if $x,u_i,u_j$ are as above we have for $s\iin G$ that
\[
(u_iu_j)\mult x(s)=\langle s\mult(u_iu_j),x\rangle
=\langle s\mult u_j\,s\mult u_i,x\rangle
=\langle s\mult u_j, \kappa^j_i(x)\rangle\langle s\mult u_i,x\rangle
=u_j\mult\kappa^j_i(x)u_i\mult x(s).
\]
Hence we find for $x_1,x_2\iin\Phi_{\ccoef(\pi_j)}$ we have for $u_i,u_j$ as above
\[
\langle u_i,\kappa^j_i(x_1x_2)\rangle
=\frac{\langle u_iu_j, x_1x_2\rangle}{\langle u_j,x_1x_2\rangle}
=\frac{\langle u_j\mult\kappa^j_i(x_1)u_i\mult x_1,x_2\rangle}{\langle u_i\mult x_1,x_2\rangle}
=\langle u_i,\kappa^j_i(x_1)\kappa^j_i(x_2)\rangle.
\]

Finally, the structure of $\Phi_{\ccoef(\pi_J)}$ and of the multiplication
on $\Phi_{\ccoef(\pi_J)}\cup\{0\}$, can now be deduced from the results
indicated in the proof of Proposition \ref{prop:ilies}.  \endpf

\subsection{Abelian groups}
For this section let $G$ denote a locally compact abelian group.  We let $\what{G}$
denote the dual group and $\bl^1(\what{G})$ its group algebra.
If $U\subset\what{G}$ is any open subset, we let $\bl^p(U)
=\{f\in \bl^p(\what{G}):f=1_Uf\}$, for $p=1,2,\infty$.  We define
\[
\pi_U:G\to\fU(\bl^2(U))\text{ by }\pi_U(s)f(\chi)=\chi(s)f(\chi)=\hat{s}(\chi)f(\chi).
\]
If $U=G$ then $\pi_G\cong\lam$ via conjugation by the Plancherel unitary,
hence $\pi_U\qc\lam$. In fact, if we let
$\lam_U=\int^\oplus_U \chi\,dm(\chi)$ on $L^2(U)=\int^\oplus_U \Cee_\chi\,dm(\chi)$
where $m$ is the Haar measure on $\what{G}$,
then $\pi_U\cong\lam_U$.
The Fourier transform gives $\ccoef(G)\cong\bl^1(G)$ and
restricts to give the identification $A_{\pi_U}\cong\bl^1(U)$.
If $S=\bigcup_{n=1}^\infty U^n$ is the semigroup generated by $U$, we have that
$\bl^1(S)$ is the closed algebra generated by $U$, hence $\ccoef(\pi_U)=\ccoef_{\pi_S}$,
and we thus consider the algebra
\begin{equation}\label{eq:piS}
\ccoef(\lam_S)=\ccoef(\pi_S)\cong\bl^1(S).
\end{equation}


We let $\Dee=\{z\in\Cee:|z|\leq 1\}$ denote
the closed disc; while this notation differs from that in most complex analysis texts,
it is more convenient for our needs.  We let for any open semigroup $S\subset\what{G}$
\[
\what{S}=\{\sig:S\to\Dee\,|\,
\sig\text{ is continuous and }\sig(\chi\chi')=\sig(\chi)\sig(\chi')\ffor \chi,\chi'\iin S\aand\sig\not=0\}
\]
denote the set of bounded semicharacters, which is itself an involutive
semigroup under pointwise multiplication and conjugation.


\begin{proposition}\label{prop:semicharacters}
If $S$ is an open subsemigroup of $\what{G}$, then $\Phi_{\ccoef(\pi_S)}
\cong \what{S}$.
\end{proposition}


\proof An elegant way to prove this result is to use the identification (\ref{eq:piS})
and prove directly that $\what{S}\cong\Phi_{\bl^1(S)}$ via the identification
$\sig\mapsto\left(f\mapsto\int_S f(\chi)\sig(\chi)\,dm(\chi)\right)$.  The later identification
is shown in \cite[4.1]{arenss}.
However, we wish to emphasise how this follows from Theorem \ref{theo:walter}.

We first observe that $\vn(\pi_S)=\{M_\vphi:\vphi\in\bl^\infty(S)\}$,
where each operator $M_\vphi$ is given by $M_\vphi f=\vphi f$.
Indeed, it is well known that $\spn\{\hat{s}:s\in G\}$ is weak* dense in
$\bl^\infty(\what{G})$, and it follows that $\spn\{1_S\hat{s}:s\in G\}$ is weak* dense in
$\bl^\infty(S)=1_S\bl^\infty(\what{G})$.    Since $\pi_S\otimes\pi_S\qc\pi_S$,
we have, in the notation of Theorem \ref{theo:walter}, a normal $*$-homomorphism
$(\pi_S\otimes\pi_S)^{\pi_S}:\vn_{\pi_s}\to\vn_{\pi_S\otimes\pi_S}
\cong\vn_{\pi_S}\wbar{\otimes}\vn_{\pi_S}$.  Identifying $\bl^2(S)\otimes^2\bl^2(S)
\cong\bl^2(S\cross S)$ we compute that
$\pi_S\otimes\pi_S(s)=M_{\hat{s}\comp \varsigma}$ where $\varsigma:S\cross S\to S$
is the multiplication map.  Thus for $\vphi\iin\spn\{1_S\hat{s}:s\in G\}$,
$(\pi_S\otimes\pi_S)^{\pi_S}(M_\vphi)=M_{\vphi\comp\varsigma}$, and by weak*
continuity this identification extends to all $\vphi\in\bl^\infty(S)$.  If $\sig\in\bl^\infty(S)$,
$M_\sig\in\Phi_{A(\pi_S)}\cup\{0\}$ if and only if $\sig$ has essential range within $\Dee$
and, by Theorem \ref{theo:walter}, $\sig\comp \varsigma
=\sig\otimes\sig$, i.e.\ $\sig(\chi\chi')=\sig(\chi)\sig(\chi')$ for a.e.\ $\chi,\chi'\iin S$.
We note that for
$f\in\bl^1(S)$, $\int_S\sig(\chi')f(\chi^{-1}\chi')\, dm(\chi')
=\sig(\chi)\int_S \sig(\chi')f(\chi')\,dm(\chi')$, from which it is
immediate that $\sig$ is continuous.
\endpf

The algebras $\ccoef(\pi_S)$ above, are all algebras of generalised
analytic functions in the sense of \cite{arenss}.

It is noted in \cite{arenss} that the decomposition
$\Dee\setdif\{0\}=(0,1]\cross\Tee$ gives a polar decomposition
$\sig(t)=|\sig(t)|\mathrm{sgn}\sig(t)$.  In light of the identification
(\ref{eq:piS}), this is a special case of the polar decomposition observed in
Theorem \ref{theo:walter}.

For concreteness, we record some particular examples of semicharacter semigroups
whose descriptions we have not been able to locate in the literature.
For subgroups of vector groups we always use additive notation.

\begin{proposition}\label{prop:someduals}
{\bf (i)}  Let $S$ be any subsemigroup of the additive semigroup $\Zee^{\geq 0}$.
Then
\[
\what{S}\cong\begin{cases} \Dee &\iif 0\in S \\
\Dee\setdif\{0\} &\iif 0\not\in S\end{cases}
\]
where each semicharacter is given for $s\iin S$
by $\sig_z(s)=z^{s/d}$ for some $z\iin\Dee$,
where $d=\mathrm{gcd}(S\setdif\{0\})$.

{\bf (ii)} Let $S$ be any open subsemigroup of the vector group $\Ree^n$.
Then there is a linearly independent
subset $\{h_1,\dots,h_n\}$ of $\Ree^N$ and $0\leq l\leq n$ for which
$S\subset\bigoplus_{j=1}^l\Ree h_j\oplus\bigoplus_{j=1}^{n-l}\Ree^{>0}h_j$.
Moreover $\what{S}\cong\Ree^l\cross\Hee^{n-l}$ where $\Hee=\{z\in\Cee:\Im z\geq 0\}$
and each semicharacter is given by
\[
\sig_{x,z}(s)=e^{i(x_1s_1+\dots+x_ls_l+z_1s_{l+1}+\dots+z_{n-l}s_n)}
\]
where $x\in\Ree^l$, $z\in\Hee^{n-l}$ and $s=\sum_{j=1}^ns_jh_j\in S$.
\end{proposition}

\proof {\bf (i)}  We may write $S\setdif\{0\}=\{s_1,s_2,\dots\}$ where $s_k<s_{k+1}$
for each $k$.  We note that $d_k=\mathrm{gcd}(s_1,\dots,s_k)$ defines
a decreasing sequence of natural numbers, hence there is $n$ such that
$d_n=d$.  The ``postage stamp problem'' tells us that there is $s_0\iin S\setdif\{0\}$
such that $S_0=\{s_0+kd:k\in\Zee^{\geq 0}\}=\{s\in S:s\geq s_0\}$,
i.e.\ if $d=m_1s_1+\dots+m_ns_n$
where $m_1,\dots,m_n\iin\Zee$ are found by the Euclidean algorithm, then
$s_0=(|m_1|+\dots+|m_n|)s_1\dots s_n$ suffices.

Now let $\sig\in\what{S}$.  We note for any $s,t\iin S\setdif\{0\}$ we have
$\sig(t)^s=\sig(st)=\sig(s)^t$, hence if $\sig(s)=0$ for any $s\iin S\setdif\{0\}$
then $\sig|_{S\setdif\{0\}}=0$; let us assume otherwise.
We define $\tau:\En\to\Cee$ by
\[
\tau(k)=\frac{\sig(s+kd)}{\sig(s)}
\]
for some $s\iin S_0$.  This definition is independent of  the choice of
such $s$ since $\sig(s+kd)\sig(t)=\sig(s)\sig(t+kd)$.  Thus we find that for
$k,l\iin\En$ that
\[
\tau(k+l)=\frac{\sig(s+kd+ld)}{\sig(s)}\frac{\sig(s+ld)}{\sig(s+ld)}=\tau(k)\tau(l)
\]
and hence $\tau(k)=\tau(1)^k$ for $k\iin\En$; let $z=\tau(1)$.  Now for
$s\iin S_0$ we write $s=kd$ for some $k\iin\Zee^{\geq 0}$ and we have
\[
\sig(s)=\frac{\sig(s+kd)}{\sig(s)}=z^k=z^{s/d}.
\]
If $s\in S\setdif (S_0\cup\{0\})$ let $k$ be so $kd=s$. Then there is $m_0\in\En$
for which $m\geq m_0$ implies that
$ms\in S_0$.  For any such $m$, $\sig(s)^m=\sig(ms)=(z^{s/d})^m$, from
which it follows that $\sig(s)=z^{s/d}$.

{\bf (ii)}  Let $S_0=\{s\in S:\Ree^{\geq 1} s\in S\}$.
It is clear that $S_0$ is a subsemigroup of $S$.  Fix a norm $|\cdot|$ on $\Ree^n$
and for any $x\iin\Ree^n$ and $\del>0$ let $B_\del(x)$ denote the open $\del$-ball
about $x$.  Since for $t\iin S$ there is $\del>0$ so $B_\del(t)\in S$,
$\bigcup_{k=1}^\infty B_{k\del}(kt)\subset S$, and we have
\begin{equation}\label{eq:mtinso}
mt\in S_0\ffor m\iin\En\text{ with }m\del \geq |t|.
\end{equation}
Moreover $S_0$ is open.  Indeed fix $s_0\iin S_0$.
By compactness, $\del=\min\{\mathrm{dist}([1,2+|s_0|]s_0,\Ree^n\setdif S),1\}>0$,
and it follows for $s$ with $|s-s_0|<\del/(1+|s_0|)$, that $[1,1+|s|]s\in S$.
Thus, by a straighforward tiling argument, $s\in\Ree^{\geq 1}s\subset S_0$.

Now let $D=\{h\in\Ree^n:s+\Ree^{\geq 0} h\in S_0\text{ for all }s\iin S_0\}$
It is clear that $D$ is a semigroup of $\Ree^n$,
$S_0\subset D$ and $\Ree^{\geq 0}h\subset D$ for $h\iin D$.
It follows that $D-D$ is a subgroup of $\Ree^n$ with non-empty interior,
and hence all of $\Ree^n$, thus $D$ contains a basis $\{h_1,\dots,h_n\}$
for $\Ree^n$.  We may re-order the basis and let $l$ be so
$\Ree h_j\subset D$ for $j=1,\dots,l$ and $\Ree h_j\not\subset D$
for $j=l+1,\dots,n$.  We then find $D=\bigoplus_{j=1}^l\Ree h_j\oplus
\bigoplus_{j=1}^{n-l}\Ree^{\geq 0} h_j$.  We let $H$ denote the interior of $D$.
Now it is well-known that $\what{\Ree^{>0}}\cong\Hee$ via the identification
$\sig(t)=e^{itz}$. 
Since the formula $\what{S_1\cross S_2}\cong\what{S}_1\cross\what{S}_2$
holds, we find that $\what{H}=\Ree^l\cross\Hee^{n-l}$ with
duality as suggested above.

Now let $\sig\in\what{S}$.   Let $\tau:H\to\Cee$ be given by
\[
\tau(h)=\frac{\sig(s+h)}{\sig(s)}
\]
for any $s\iin S_0$ for which $\sig(s)\not=0$; by (\ref{eq:mtinso})
such an $s$ always exists.  As in the proof
of (i), above, $\tau\in\what{H}$ and is independant of choice of $s$.
We note that $S\in H$.
Indeed, if $t\in S$ pick $m$ as in (\ref{eq:mtinso}).  If $s\in S_0$
and $\alp\geq 1$ we have $\alp(s+t)=\alp s+\frac{\alp}{m}mt$
and we note that $S_0\in H$.  It is immediate that
for $t\in S$ that $\tau(t)=\sig(s+t)/\sig(s)=\sig(t)$.  \endpf

\begin{example} \label{ex:algx}
Let $G=\Tee$ and $\chi_n\iin\what{\Tee}$ be given by $\chi_n(z)=z^n$
for $z\iin\Tee$.   For any subset $U$ of $\Zee$,
let $\lam_U=\bigoplus_{n\in U}\chi_s$ where $\chi_s(z)=z^s$.  Then
\[
\ccoef_{\lam_U}=\left\{z\mapsto\sum_{n\in U}\alp_n z^n\text{ where }
\sum_{n\in U}|\alp_n|<\infty\right\}\cong\ell^1(U)
\]
is a Banach space of Laurent polynomials. Let $S$
be the subsemigroup generated by $U$.  We have essentially two cases to consider.

{\bf (i)} If $S$ is a semigroup which contains both positive
and negative elements, it is a group and hence of the form $\Zee d$.
In this case $\Phi_{\ccoef(\lam_S)}\cong\Tee$, via the character
$\chi(z)=z^d$.

{\bf (ii)} If $U=\{0,1\}$ then $\tau_{\lam_U}=\lam_{\Zee^{\geq 0}}$, and we find that
$\fA(\lam_S)=\fA_{\lam_{\Zee^{\geq 0}}}=\fA(\Dee)$ is the disc algebra,
consisting of functions on $\Tee$ which are continuous and continuously extend
to analytic functions on the interior of $\Dee$.
If $S$ is a subsemigroup of $\Zee$ which is not a group, then we
may suppose that $S\subset\Zee^{\geq 0}$, otherwise take $-S$.
Proposition \ref{prop:someduals} (i) shows that
if $U$ is any subset of $\En$, then $\fA(\lam_S)$ may be identified with
the uniformly closed subalgebra of analytic functions on $\Dee$ generated by
the monomials $z\mapsto z^n$ for $n\in U$.

{\bf (iii)} Note that if we consider the topological semigroup $\Dee$ itself,
it is an obvious consequence of the maximum modulus principle that
the translation invariant algebra $\fA(\Dee)$ has \v{S}ilov boundary
$\Tee$, which is not an ideal in $\Dee$.
Thus Theorem \ref{theo:algspec} (iii) may not be true for any
semitopological semigroup, without assuming the existence of a
dense subgroup.
\end{example}

\begin{example}\label{ex:realline}
Let $G=\Ree$ and $\chi_s\iin\what{\Ree}$ be given by $\chi_s(t)=e^{ist}$.

{\bf (i)} Let $\lam_+=\int^\oplus_{\Ree^{>0}}\chi_s\, dm(s)$ where $m$ is Lebesgue measure.
Then it is standard that
$\ccoef_{\lam_+}=\ccoef(\lam_+)\cong \bl^1(\Ree^{>0})$.  Propositions
\ref{prop:semicharacters} and \ref{prop:someduals} (ii) show that the Gelfand
transform on $\bl^1(\Ree^{>0})$ is given by the Laplace transform
\[
\hat{f}(z)=\int_0^\infty f(t)e^{izt}dt
\]
where $z\in\Hee$; it is clear that $\hat{f}$ is analytic on the interior of $\Hee$
and vanishes at $\infty$.  Let $\fA_0(\Hee)$ denote the algebra of continuous
functions on $\Ree$, which continuously extend to analytic functions
on the interior of $\Hee$ which vanish at $\infty$ on all of $\Hee$.
Note that $\Phi_{\fA_0(\Hee)}\cong\Hee$ and $\partial_{\fA_0(\Hee)}=\Ree$.
Let us show that the uniform closure of
$\{\hat{f}|_\Ree:f\in\bl^1(\Ree^{>0})\}$ is $\fA_0(\Hee)$.
First we consider the Cayley transform $\gam:\Hee\to\Dee\setdif\{1\}$
given by $\gam(z)=\frac{z-i}{z+i}$.  The map
$g\mapsto g\comp\gam:(1-z)\fA(\Dee)\to\fA_0(\Hee)$ is an isomorphism.
If $g_n(z)=z^n-z^{n+1}$ then $g_n\comp\gam(z)=\frac{2(z-i)^n}{i(z+i)^{n+1}}$.
Let $f_n\in\bl^1(\Ree^{>0})$ be given by $f_n(t)=t^ne^{-t}$.  Then
$\hat{f}_n(z)=\int_0^\infty t^ne^{i(z+i)t}dt=-\frac{i^nn!}{(z+i)^{n+1}}$.
Hence it follows that $\spn\{\hat{f}_k\}_{k=1}^n=\spn\{g_k\comp\gam\}_{k=1}^n$
for each $n\iin\En$, so $\{\hat{f}:f\in\bl^1(\Ree^{>0})\}$ is dense
in $\fA_0(\Hee)$.

{\bf (ii)}  If $S$ is any open subsemigroup of $\Ree^{>0}$ we let
$\lam_S=\int^\oplus_S\chi_s\,dm(s)$.  As above we get $\ccoef_{\lam_S}=
\ccoef(\lam_S)\cong\bl^1(S)$.  The Gelfand transform on $\bl^1(S)$,
in this case, is a modified Laplace transform
\[
\hat{f}(z)=\int_S f(t)e^{izt}dt.
\]
By Proposition \ref{prop:someduals} (ii) and Theorem \ref{theo:algspec} (i),
$\fA(\lam_S)$ is isomorphic to a $\Hee$-translation invariant subalgebra of
$\fA_0(\Hee)$.

The simplest example is the semigroup $\Ree^{>a}$ where
$a>0$.  By standard Laplace transform techniques we see that
$\hat{f}(z)=e^{iaz}\int_0^\infty f(t-a)e^{izt}dt$ from which it follows
that $\fA(\lam_{\Ree^{>a}})\cong e^{iaz}\fA_0(\Hee)$.

If we let $\frac{1}{5}<a<\frac{1}{3}$, then the semigroup
$S=\{s\in\Ree:1-a<a<1+a\text{ or }s>2-2a\}$ can be shown, as above,
to satisfy $\fA(\lam_S)=(e^{i(1-a)z}-e^{i(1+a)z}+e^{iaz})\fA_0(\Hee)$.

Notice that in both cases above $\fA(\lam_S)$ is a principal
ideal in $\fA(\lam_{\Ree^{>0}})$.

{\bf (iii)} In \cite{grigoryant} many examples of the form
$\pi=\bigoplus_{s\in S}\chi_s$, where $S$ is a subsemigroup of
$\Ree^{\geq 0}$, are given.  These correspond to certain analytic
semigroups which contain quotients of $\Ree^{\fA\fP}$.

{\bf (iv)} Consider the representation $\pi_+\oplus \chi_1$.  We have
that $\tau_{\pi_+\oplus \chi_1}\cong_q\pi_+\oplus\bigoplus_{n\in\En}\chi_n$
so
\[
\ccoef(\lam_+\oplus \chi_1)=\ccoef_{\lam_+}\oplus_{\ell^1}\ccoef(\chi_1)
\cong\bl^1(\Ree^{>0})\oplus_{\ell^1}\ell^1(\En)
\]
where $\ell^1(\En)$ is the algebra of Dirac measures supported on $\En$.
Taking uniform closure, we obtain a semidirect product algebra
\[
\fA(\lam_+\oplus \chi_1)=\fA_0(\Hee)\oplus\fA_0(\Dee)
=\fA_0(\Hee\sqcup\Dee_0)
\]
where $\fA_0(\Dee)=z\fA(\Dee)$.
As suggested by Proposition \ref{prop:morespine},
 $\Hee\sqcup\Dee_0$ is a semigroup where $\Hee$
and $\Dee_0=\Dee\setdif\{0\}$ are subsemigroups, and
for $z\iin\Hee$ and $w\iin\Dee_0$ we define
\[
zw=e^{iz}w=wz.
\]
The topology is given by having $\Hee$ be open, and allowing neighbourhoods
of elements $w\iin\Dee_0$ to be given by $U\sqcup V$ where $V$ is a neighbourhood
of $w$ in $\Dee_0$, and $U$ is a cocompact set in $\Hee$.

{\bf (v)} Let $S=\Ree\cross\Ree^{>0}\subset\Ree^2$.  This is a
subgroup for which we have
\[
\fA(\lam_S)\cong\fC_0(\Ree)\check{\otimes}\fA_0(\Hee)
\]
where $\check{\otimes}$ denotes the injective tensor product.
Indeed, we have an isometric identification
$\fC_0(\Ree)\check{\otimes}\fA_0(\Hee)\cong\{u\in\fC_0(\Ree^2):
u(x,\cdot)\in\fA_0(\Hee)\}$.  Since $\ccoef(\lam_S)\cong\bl^1(\Ree\cross\Ree^{>0})
\cong\bl^1(\Ree)\hat{\otimes}\bl^1(\Ree^{>0})$ (projective tensor product),
it follows, in part from (i) above,
that for each $f\iin\bl^1(S)$,  the ``Laplace-Fourier transform" of $f$
\[
\hat{f}(x,z)=\int_S f(s)e^{i(s_1x+s_2z)}\,ds
\]
satisfies $f|_{\Ree^2}\in\fC_0(\Ree)\check{\otimes}\fA_0(\Hee)$,
and the family of such functions is uniformly dense within.

Similarly, if $a_1,\dots,a_{n-l}\geq 0$ and
$S=\Ree^l\cross\Ree^{>a_1}\cross\Ree^{>a_{n-l}}\subset\Ree^n$,
then we can appeal to (ii) above and adapt the above methods to see that
\[
\fA(\lam_S)\cong\fC_0(\Ree^k)\check{\otimes}e^{ia_1z}\fA_0(\Hee)
\check{\otimes}\dots\check{\otimes}e^{ia_{n-l}z}\fA_0(\Hee).
\]
\end{example}

\begin{example}\label{ex:elgun}
We note some recent results of \cite{elgun}.  For $G=\Zee$, there is a representation
$\pi$ for which $G^\pi\cong\bl^\infty[0,1]_{\norm{\cdot}\leq 1}$.

Furthermore, if $G=(\Zee/p\Zee)^{\oplus\infty}$, then
there is a representation $\pi$ for which $G^\pi
\cong\{\vphi\in\bl^\infty[0,1]:\mathrm{ess\,ran}\vphi\subset \Gamma_p\}$
where $\Gamma_p=\mathrm{conv}\{e^{2\pi ik/p}:k=0,\dots,p-1\}$.

More significantly, for each example above
there are $\mathfrak{c}$ idempotents
in $G^\pi$, and it is shown the closure of this family is isomorphic to
$\{\vphi\in\bl^\infty[0,1]:\mathrm{ess\,ran}\vphi\subset[0,1]\}$.
\end{example}

\subsection{Compact matrix groups}

Let $G$ be a compact group, and $\sig$ be a continuous finite dimensional unitary
representation, so $G^\sig=\sig(G)$ may be regarded as a closed
subgroup of the unitary group $\fU(\fH_\sig)\subset\fB(\fH_\sig)$.

We now use some ideas form Lie theory.  We let
\[
\gg^\sig=\{X\in \fB(\fH_\sig):\exp(tX)\in G^\sig\text{ for all }t\iin\Ree\}.
\]
It is well-known that $\gg^\sig$ is a real Lie algebra with $[X,Y]=XY-YX$.  Moreover
$\gg^\sig$ has the same reducing subspaces as $\sig$, so $\gg^\sig\subset\vn_\sig$.
We let $\gg^\sig_\Cee=\gg^\sig+i\gg^\sig$ denote its complexification and
\[
G^\sig_\Cee=G^\sig\langle \exp\gg^\sig_\Cee\rangle 
\]
which is a complex Lie subgroup of invertible elements in $\vn_\sig$.  We then let
$\Dee^\sig=\wbar{G^\sig_\Cee}\cap\fB(\fH_\sig)_{\norm{\cdot}\leq 1}$

\begin{theorem}\label{theo:sigspectrum}
We have $\Phi_{\ccoef(\sig)}\cong \Dee^\sig\setminus\{0\}$.
\end{theorem}

\proof For this proof, we shall make use of the Zariski topolgy
on the finite dimensional affine space $\fB(\fH_\sig)$.  We let
$\pol(\fB(\fH_\sig))$ denote the algebra generated by the
matrix coefficient functionals and the constant functional $1$.  Then
for $S\subset\fB(\fH_\sig)$, we let $\ideal(S)=\{p\in\pol(S):p|_S=0\}$
and let the Zariski closure of $S$ be given by
\[
Z(S)=\bigcap_{p\in\ideal(S)}p^{-1}(\{0\}).
\]
We observe that $Z(\{0\})=\{0\}$ and  $Z(S\cup\{0\})=Z(S)\cup\{0\}$
since $S\mapsto Z(S)$ is a closure operation.  We also note that
$Z(S)$ is the largest set $Z$ for which $\ideal(S)=\ideal(Z)$ and hence
the spectrum of the algebra $\pol(\fB(\fH_\sig))|_S\cong
\pol(\fB(\fH_\sig))/\ideal(S)$ is naturally identified with $Z(S)\setminus\{0\}$.
Thus, recognising $\alg(F_\sig)$ as the algebra $\pol(\fB(\fH_\sig))|_{G^\sig\cup\{0\}}$,
we obtain spectrum 
\begin{equation}\label{eq:specisz}
\Phi_{\alg(F_\sig)}\cong Z(G^\sig\cup\{0\})\setminus\{0\}=Z(G^\sig)\setminus\{0\}.
\end{equation}

We next wish to establish that $G^\sig_\Cee\subset Z(G^\sig)$.
First, if $p\in\ideal(G^\sig)$, we have that
$p\comp\exp:\gg^\sig_\Cee\to\Cee$ is a holomorphic function, which vanishes
on $\gg^\sig$, a real subspace whose complex span is $\gg^\sig$. 
Hence $\exp(\gg^\sig_\Cee)\subset Z(G^\sig)$. 
By virtue of the fact that $G^\sig$ is
a (semi)group, we have for $p$ in $\ideal(G^\sig)$ that 
$s\mult p\in \ideal(G^\sig)$ for $s$ in $G^\sig$.  Hence for $v$ in $Z(G^\sig)$,
$p(vs)=0$ for $p$ and $s$ as above, and we have that $p\mult v\in\ideal(G^\sig)$.  Thus
for $w$ in $Z(G^\sig)$ we find $p(vw)=p\mult v(w)=0$, so $vw\in Z(G^\sig)$.
Thus $Z(G^\sig)$ is a semigroup, and it follows that $G^\sig_\Cee=
G^\sig\langle\exp(\gg^\sig_\Cee)\rangle\subset Z(G^\sig)$.

We now wish to establish that $G^\sig_\Cee$ is Zariski closed in $\fB(\fH_\sig)_{inv}$.
To see this we consider $G^{\sig\oplus\bar{\sig}}_\Cee$.  We observe that
\[
G^{\sig\oplus\bar{\sig}}\cong\left\{\begin{bmatrix} u & 0 \\ 0 & \bar{u}\end{bmatrix}
:u\in G^\sig\right\}\cong G^\sig
\]
Hence, following calculations such as in \cite[Cor.\ 2.2]{ludwigst}, we see that
\[
G^{\sig\oplus\bar{\sig}}_\Cee\cong\left\{\begin{bmatrix} v & 0 \\ 0 & v^{-T}
\end{bmatrix}:v\in G^\sig_\Cee\right\}\cong G^\sig_\Cee
\]
where $v^{-T}$ is the inverse transpose of $v$.  We also note that we can identify
$\alg(F_{\sig\oplus\bar{\sig}})$ with the algebra of trigonometric functions on $G^\sig$,
thanks to \cite[(27.39)]{hewittrII}.  Then $G^\sig_\Cee$ is naturally identified
with the spectrum of $\alg(F_{\sig\oplus\bar{\sig}})$, \cite[(2.3) \& Thm.\ 2]{cartwrightm}.
Combining this with (\ref{eq:specisz}) we obtain that $G^{\sig\oplus\bar{\sig}}_\Cee=
V(G^{\sig\oplus\bar{\sig}})$.   The same polynomial equations, on matrices in the 
upper left corner, establish that $G^\sig_\Cee$ is Zariski closed in $\fB(\fH_\sig)_{inv}$.

It follows from the fact above that $
G^\sig_\Cee=Z(G^\sig_\Cee)\cap\fB(\fH_\sig)_{inv}=Z(G^\sig)\cap\fB(\fH_\sig)_{inv}.$
However, since $\fB(\fH_\sig)_{inv}$ is Zariski open, we can use
\cite[I.10 Thm.\ 1]{mumford}, to establish that
\[
\wbar{G^\sig_\Cee}=\wbar{Z(G^\sig)\cap\fB(\fH_\sig)_{inv}}=Z(G^\sig).
\]

Finally, we note that each element of $\Phi_{\ccoef(\sig)}$ is 
contractive and determined by its restriction to the dense
subalgebra $\alg(F_\sig)$; while each contractive element of 
$\wbar{G^\sig_\Cee}\cong\Phi_{\alg(F_\sig)}$ extends to a character
on $\ccoef(\sig)$.  \endpf

We note that by \cite[Cor.\ 1]{mckennon}, each $v$ in
$G^\sig_\Cee$ admits a polar decomposition $u|v|$ where $u\in G^\sig$
and hence $v\in G^\sig_\Cee$.
This corresponds to the polar decomposition observed in Theorem \ref{theo:walter}.

\begin{example}\label{ex:unitary}
{\bf (i)} Let $G=\U(d)$, the group of $d\cross d$-unitary matrices, and let $\sig:\U(d)\to
\fU(\Cee^d)$ denote the standard representation, so $G^\sig=\fU(\Cee^d)$.  It is well-known
that the Lie algebra is $\mathfrak{u}(d)=\{X\in\fB(\Cee^d):X^*=-X\}$, whence
$\mathfrak{u}(d)_\Cee=\fB(\Cee^d)$, and hence $\fU(\Cee^d)_\Cee=\fB(\Cee^d)_{inv}$.
This space is dense in $\fB(\Cee^d)$, so by Theorem \ref{theo:sigspectrum}
$\Phi_{\ccoef(\sig)}=\fB(\Cee^d)_{\norm{\cdot}\leq 1}\setminus\{0\}$.

It is standard that the convex hull of $\fU(\Cee^d)$ is
$\fB(\Cee^d)_{\norm{\cdot}\leq 1}$, i.e.\ we have for $v$ in $\fB(\Cee^d)_{\norm{\cdot}\leq 1}$
polar decomposition $v=u|v|=\frac{1}{2}u[(|v|+i\sqrt{1-|v|^2})+(|v|-i\sqrt{1-|v|^2})]$.
Hence, each element of $\fB(\Cee^d)_{\norm{\cdot}\leq 1}$ may be veiwed as
a convex combination of elements $\sig^\sig_e(\eps_{\fA(\sig)}(u))$ --- see notation of
Corollary \ref{cor:walter1} --- for $u$ in $\U(d)$.   It follows
that $\Phi_\fA(\sig)\cong\fB(\Cee^d)_{\norm{\cdot}\leq 1}$, thus verifying
Conjecture \ref{conj:nonwp} in this case.  We obseve that
$\fA(\sig)\cong\fA_0(\fB(\Cee^d)_{\norm{\cdot}\leq 1})$, the ``punctured ball algebra", i.e.
the algebra of all
continuous functions on $\fB(\Cee^d)_{\norm{\cdot}\leq 1}$ which are holomorphic
on the interior and vanish at $0$.  Similarly we find that
$\fA(\sig\oplus 1)\cong\fA(\fB(\Cee^d)_{\norm{\cdot}\leq 1})$, the ball algebra.
It is well-known that the \v{S}ilov boundary of $\fA(\fB(\Cee^d)_{\norm{\cdot}\leq 1})$
is $\fU(\Cee^d)$, see, for example \cite[\S 12.4]{faraut}.

{\bf (ii)} Let $H$ denote either of the classical compact matrix groups $\O(d)$
or $\Sp(d)$.
Let $G=\Tee\mult H=\{zu:z\in\Tee\text{ and }u\in H\}$, and
$\sig:G\to\fB(\Cee^{d'})$ be the standard representation where $d'=d$
in the case that $H=\O(d)$, and $d'=2d$ in the case that $H=\Sp(d)$.  
It is straightforward to compute that
$G^\sig_\Cee=\Cee^{\not=0}\mult H_\Cee=\{\alp v:\alp\in\Cee^{\not=0}
\text{ and }v\in H_\Cee\}$.  We note that
$\O(d)_\Cee=\O(d,\Cee)$ and $\Sp(d)_\Cee=\Sp(d,\Cee)$ are each closed algebraic
groups.

We claim that $\wbar{G^\sig_\Cee}=\Cee\mult H_\Cee$.  
Indeed, if $\lim_{n\to\infty}\alp_n v_n=b$ then we factor
$\alp_n v_n=(\alp_n\norm{v_n})(\frac{1}{\norm{v_n}}v_n)$.
By dropping to a subsequence, we may suppose that
$\lim_{n\to\infty}\alp_n\norm{v_n}=\alp$ and
$\lim_{n\to\infty}\frac{1}{\norm{v_n}}v_n=v\in H_\Cee$.  Hence $b=\alp v\in
\Cee\cdot H_\Cee$.  Thus $\Dee^\sig=\Cee\mult H_\Cee\cap \fB(\Cee^{d'})_{\norm{\cdot}\leq 1}$.

Thus $\Phi_{\ccoef(\sig)}\cong \Dee^\sig\setminus\{0\}
=G^\sig_\Cee\cap\fB(\Cee^{d'})_{\norm{\cdot}\leq 1}$.
We have not devised a means to show that $\Phi_{\fA(\sig)}=\Phi_{\ccoef(\sig)}$
in either of these cases.


{\bf (iii)} If $\sig:G\to\fU(\fH_\sig)$ is an injective homomorphism, then
$\rho_\sig\cong_q\lam$, the left regular representation, by \cite[(27.39)]{hewittrII}
and the Peter-Weyl Theorem.  Thus $\ecoef(\sig)=\falg$ and
$\Phi_{\ecoef(\sig)}=G$.  

In particular if $G=\mathrm{SU}(2)$ and $\sig:\mathrm{SU}(2)\to\U(2)$
is the standard representation, then $\bar{\sig}\cong\sig$ and it follows that
$\ccoef(\sig)=\ecoef(\sig)=\ccoef(\mathrm{SU}(2))$.

\end{example}

\subsection{A non-compact, non-abelian example}

\begin{example}\label{ex:axplusb}
Let $G$ be the $ax+b$-group, given by 
\[
\{(a,b)
:a\in\Ree^{>0},b\in\Ree\}
\]
with multiplication $(a,b)(a',b')=(aa',ab'+b)$.
Let us consider, adapted from the notation of
\cite[p.\ 189]{folland} (we omit normalisation by $2\pi$), the
representation $\pi_+:G\to\fB(\bl^2(\Ree^{>0},m))$ given by
\begin{equation}\label{eq:piplus}
\pi_+(a,b)f(s)=a^{1/2}e^{ibs}f(as)
\end{equation}
for $f\iin\bl^2(\Ree^{>0},m)$ and $m$-a.e.\ $s\iin\Ree^{>0}$, where $m$ is the usual
Lebesgue measure.     We wish to compute $\Phi_{\fA(\pi_+)}$.

According to \cite[Th\'{e}o.\ 5]{khalil}, 
$\ccoef_{\pi_+}=\ccoef(\pi_+)$.  Implicit in the proof of that fact, see \cite[p.\ 159]{khalil}
is a formula for $\pi_+\otimes\pi_+$.  
We rederive this formula in a form more tractable to to obtaining  (\ref{eq:axbwalter}).
For convenience
we write $\bl^2=\bl^2(\Ree^{>0},m)$, which we identify as the subspace
$\bl^2(\Ree)1_{\Ree^{>0}}$ of $\bl^2(\Ree)$.  We consider the direct integral
Hilbert space $\fH=\int^\oplus_{\Ree^{>0}} \bl^2_t\,dt$ where each
$\bl^2_t$ is a copy of $\bl^2$.  
We identify $\bl^2\otimes^2\bl^2$
with $\bl^2(\Ree^2)1_{(\Ree^{>0})^2}$ in the usual manner, and define
$U:\bl^2\otimes^2\bl^2\to\fH$ by
\[
U\xi=\int^\oplus_{\Ree^{>0}}(\lam(t)\otimes I)\xi(\cdot,t)\,dt
\]
where $\lam(t)f(s)=f(t-s)$ for $f\in\bl^2(\Ree)1_{\Ree^>0}\cong\bl^2(\Ree^{>0})$.
It is straightforward to verify that $U$ is a unitary with
$U^*\left(\int^\oplus_{\Ree^{>0}}u_{t'}\,dt'\right)(s,t)=u_t(s+t)$.
We define for $(a,b)\iin G$, $\Pi_+(a,b):\fH\to\fH$ by
\[
\Pi_+(a,b)\int^\oplus_{\Ree^{>0}}u_{t}\,dt
=a^{1/2}\int^\oplus_{\Ree^{>0}}\pi_+(a,b)u_{at}\,dt.
\]
We have for $f,g,h,k\iin\bl^2$ that
\begin{align}
\langle \pi_+\otimes\pi_+(a,b)f\otimes g &| h\otimes k\rangle
=\inprod{\pi_+(a,b)f}{h}\inprod{\pi_+(a,b)g}{k} \notag \\
&=a\int_{\Ree^{>0}} \int_{\Ree^{>0}} e^{ib(s+t)} f(as)g(at)\wbar{h(s)k(t)}\,dt\,ds \notag \\
&=a^{1/2}\int_{\Ree^{>0}} \int_{\Ree^{>0}}a^{1/2}e^{ibs}
f(a(s-t))g(at)\wbar{h(s-t)k(t)}\,ds\,dt \label{eq:piplusprod} \\
&=a^{1/2}\int_{\Ree^{>0}}\inprod{\pi_+(a,b)\lam(at)f\,g(at)}{\lam(t)h\, k(t)}\,dt \notag \\
&=\inprod{\Pi_+(a,b)Uf\otimes g}{Uh\otimes k}. \notag
\end{align}
In particular $\Pi_+$ is a representation, unitarily equivalent to $\pi_+\otimes\pi_+$.
Thus Theorem \ref{theo:walter} (\ref{eq:walter}) tells us that for $x\iin\vn_\pi$
we have that $x\in\Phi_{\ccoef(\pi_+)}\cup\{0\}$ if and only if
\begin{equation}\label{eq:axbwalter}
U^*\Pi_+^{\pi_+}(x)U=x\otimes x.
\end{equation}
We remark that since $\pi_+$ is irreducible, $\vn_{\pi_+}=\fB(\bl^2)$.
We let 
\[
\til{G}=\{(a,z)=
:a\in\Ree^{>0},z\in\Hee\}
\]
which is a semigroup via $(a,z)(a',z')=(aa',az'+z)$.  
For $(a,z)\iin\til{G}$, define
$\til{\pi}_+(a,z)$ exactly as in (\ref{eq:piplus}).
A straightforward repeat of (\ref{eq:piplusprod})
shows that each $\til{\pi}_+(a,z)$ satsfies (\ref{eq:axbwalter}).
We conjecture that $\Phi_{\ccoef(\pi_+)}\cong\til{G}$.
We observe that the conjugation $(a,z)^*=(a^{-1},-a^{-1}\bar{z})$
satisfies $ \til{\pi}_+((a,z)^*)=\til{\pi}_+(a,z)^*$.  We have that
$(a,z)^*(a,z)=(1,2a^{-1}i\Im z)$.  It is clear that $\til{\pi}_+(1,2a^{-1}i\Im z)$
is the operator of multiplication by $x\mapsto e^{-2a^{-1}\Im z\,x}$ on $\bl^2$,
whose positive square root is the operator of multiplication by $x\mapsto e^{-a^{-1}\Im z\,x}$,
i.e.\ $|\til{\pi}_+(a,z)|=\til{\pi}_+(1,a^{-1}i\Im z)$.
Thus we obtain a formula for the
polar decomposition from Theorem \ref{theo:walter}:  $(a,z)=(a,\Re z)(1,a^{-1}i\Im z)$.

We now claim that
\begin{equation}\label{eq:piplusalg}
\fA(\pi_+)\cong\fC_0(\Ree^{>0})\check{\otimes}\fA_0(\Hee)
\end{equation}
where $\check{\otimes}$ denotes the injective tensor product and $\fA_0(\Hee)$
is defined in Example \ref{ex:realline} (i).
Indeed, $\fA(\pi_+)\subset\fC_0(G)\cong\fC_0(\Ree^{>0}\cross\Ree)$
is a point-separating subalgebra which has the span of functions
$(a,b)\mapsto a^{1/2}\int_{\Ree^{>0}}e^{ibs}f(as)\wbar{g(s)}\,ds$ where $f,g\in\bl^2$.
For fixed $b$, such a function is easily seen to be a generic element of $\ccoef(\Ree^{>0})$,
the space of which is dense in $\fC_0(\Ree^{>0})$.
For a fixed $a$, such a function may be seen
to be the Laplace transform $b\mapsto\hat{h}(b)$, of a generic element
$h\iin\bl^1(\Ree^{>0},m)$, a set of elements which is dense in
$\fA_0(\Hee)$, as demonstrated in Example \ref{ex:realline} (i).

We note that it is immediate that $\Phi_{\fA(\pi_+)}\cong\til{G}$.
Conjecture \ref{conj:nonwp} if true, would now tell us that
$\Phi_{\ccoef(\pi_+)}\cong\til{G}$.


\end{example}





\noindent{\bf Acknowledgements.}  The authors are grateful to the referee
for carefully reading the article and suggesting many corrections.  In particular,
we are grateful for the suggestion of using the Zariski closure in 
Theorem \ref{theo:sigspectrum}, which allowed us to upgrade this result
from a conjecture.  The first named author is also grateful to his colleagues
Ruxandra Moraru, David McKinnon and Rahim Moosa for valuable
discussions on the Zariski topology; and to former student
Faisal al-Faisal for a suggestion which aided an earlier version of Example
\ref{ex:unitary}

Addresses:
\linebreak
 {\sc 
Department of Pure Mathematics, University of Waterloo,
Waterloo, Ontario, N2L 3G1, Canada \\
Department of  Mathematics and Statistics,  University of Winnipeg, 
Winnipeg, Manitoba,  R3B 2E9, Canada}

\medskip
Email-adresses:
\linebreak
{\tt nspronk@uwaterloo.ca} \\ {\tt r.stokke@uwinnipeg.ca}

\end{document}